\documentclass{amsart}
\usepackage{graphicx,euscript}
\usepackage{amssymb}

\theoremstyle{plain}
\newtheorem{Thm}{Theorem}[section]
\newtheorem{Cor}[Thm]{Corollary}
\newtheorem{Lem}[Thm]{Lemma}
\newtheorem{Prop}[Thm]{Proposition}

\theoremstyle{definition}

\newcommand{\cond}{\,\lv\,}

\newcommand{\wt}{\widetilde}

\newcommand{\tp}{\tau_{\partial}}

\newcommand{\ltf}{\lim_{t\to\infty}}
\newcommand{\litf}{\liminf_{t\to\infty}}
\newcommand{\lstf}{\limsup_{t\to\infty}}
\newcommand{\lnf}{\lim_{n\to\infty}}
\newcommand{\linf}{\liminf_{n\to\infty}}
\newcommand{\lsnf}{\limsup_{n\to\infty}}
\newcommand{\lzf}{\lim_{z\to\infty}}

\newcommand{\E}{\mathbb{E}}

\newcommand{\ls}{\bigl\{}
\newcommand{\rs}{\bigr\}}
\newcommand{\te}{\exists}
\newcommand{\fa}{\forall}

\newcommand{\lp}{\bigl(}
\newcommand{\rp}{\bigr)}
\newcommand{\gm}{\gamma}
\newcommand{\gs}{\sigma}
\newcommand{\ga}{\alpha}

\newcommand{\gd}{\delta}
\newcommand{\gl}{\lambda}
\newcommand{\oh}{\frac{1}{2}}

\newcommand{\tif}{\text{ if }}

\newcommand{\gep}{\epsilon}

\newcommand{\nc}{\newcommand}

\newcommand{\gL}{\Lambda}
\newcommand{\eu}{\EuScript}
\newcommand{\indic}{\boldsymbol{1}}

\newcommand{\on}{\operatorname}
\nc{\G}{\eu{G}}
\nc{\lip}{\on{Lip}}
\nc{\izf}{\int_0^\infty}
\nc{\tand}{\text{ and }}
\nc{\tst}{\text{ s.t. }}
\renewcommand{\Pr}{\mathbb{P} \ls}
\renewcommand{\P}{\mathbb{P}}

\newcommand{\lb}{\bigl[}
\newcommand{\rb}{\bigr]}
\newcommand{\lv}{\bigl|}
\newcommand{\rv}{\bigr|}
\newcommand{\rec}{\frac{1}}

\newcommand{\R}{\mathbb{R}}
\newcommand{\Z}{\mathbb{Z}}
\nc{\N}{\mathbb{N}}
\nc{\mL}{\eu{L}}
\nc{\C}{\eu{C}}
\nc{\B}{\eu{B}}

\newcommand{\mD}{\eu{D}}
\newcommand{\mC}{\eu{C}}
\newcommand{\mc}{\mM_{1}^{c}(\R^{+})}
\nc{\mo}{\mM_{\le 1}(\R^{+})}
\nc{\supp}{\on{Supp}}
\newcommand{\mf}{\mathfrak}
\newcommand{\mI}{\mathcal{I}}

\newcommand{\sng}{\renewcommand{\baselinestretch}{1}\normalsize}

\nc{\vx}{\vec{x}}
\nc{\vy}{\vec{y}}
\nc{\DF}{\eu{F}}
\nc{\df}{f}
\nc{\tX}{\wt{X}}
\nc{\plm}{\phi_{\ulm}}
\nc{\mM}{\mf{M}}
\nc{\otn}{\phantom{}_{n}\omega_{t}}

\begin{document}
\title[Quasistationary distributions]
{Quasistationary distributions for one-dimensional diffusions with killing}
\author{David Steinsaltz and Steven N. Evans}
\address{David Steinsaltz, Department of
  Demography, University of California, 2232 Piedmont Ave., Berkeley, CA 94720\\
  Steven N. Evans, Department of Statistics, Evans Hall 367, University of
  California, Berkeley, CA 94720}
\thanks{DS supported by Grant K12-AG00981 from the National
Institute on Aging.  SNE supported in part by Grants
DMS-00-71468 and DMS-04-05778 from the National Science Foundation, and by the Miller Institute for Basic Research in Science.}

\email{dstein@demog.berkeley.edu\\evans@stat.berkeley.edu}
\date{\today}
\begin{abstract}
We extend some results on the convergence of one-dimensional 
diffusions killed at the boundary, conditioned on extended survival, to the case of general killing on 
the interior.   We show, under fairly general conditions, that a diffusion conditioned on long survival either runs off to infinity almost surely, or almost surely converges to a quasistationary distribution given by the lowest eigenfunction of the generator. In the absence of internal killing, only a 
sufficiently strong inward drift can keep the process
close to the origin, to allow convergence in distribution.  An 
alternative, that arises when general killing is allowed, is that the 
conditioned process is held near the origin by a high rate of killing 
near $\infty$.  We also extend, to the case of general killing, the standard result on convergence to a quasistationary distribution of a diffusion on a compact interval.
\end{abstract}

\maketitle

\renewcommand{\thefootnote}{*}
\section{Motivation} 
\label{sec:intro}
Among the puzzles that occupy mathematical demographers, one of the most tantalizing is the phenomenon sometimes called ``mortality
plateaus''.  While the knowledge that human mortality rates increase
with age during the years of maturity remounts to the primitive past
of statistical science, only in recent years has it become
apparent that this accelerating decrepitude slows in extreme old age,
and may even stop \cite{jV97}. Analogous flattening of the mortality
curves for fruit flies \cite{jC92} is well established.

In their attempts to explain this widespread, and possibly
near-universal phenomenon, mathematical demographers have turned
repeatedly to Markov models of mortality.  These are stochastic 
processes where ``death'' is identified with a random stopping time, 
which arises from a typically unobserved Markov process.     Some models, such as the ``cascading failures'' model
of H. Le Bras \cite{hLB76} (described in section \ref{sec:example}), were originally introduced to model the
classical exponentially increasing ``Gompertz'' mortality curve, and
were only later shown \cite{GG91} to converge to a constant plateau
mortality rate.  Others, such as the series-parallel model of L.
Gavrilov and N. Gavrilova \cite{GG91}, and the drifting Brownian
motion model of J. Weitz and H. Fraser \cite{WF01}, were introduced
explicitly with mortality plateaus in mind.  (The drifting Brownian
motion, though, it should be pointed out, was first brought into
demography by J. Anderson \cite{jA00} with slightly different motives. 
The later work of Weitz and Fraser was apparently independent.)

While these models illustrate theories of the development of 
senescence, they are not really explanations to be set beside 
the usual demographic explanations for mortality plateaus.  The 
authors derive the convergence to constant mortality from explicit 
computations alone, which cannot help to make palpable the driving force behind the convergence.  Nor 
do they suggest whether the convergence is a 
typical phenomenon, or whether it depends on the 
many arbitrary particularities of the models.

As we have explained in \cite{diffmort}, mortality plateaus in Markov 
aging models are examples of the generic property of convergence to 
quasistationary distributions.  What is more, translating 
quasistationary distributions into demographic reality offers a novel 
account of the mortality plateaus, distinct from the standard dyad of 
population heterogeneity (mortality rates stop rising because the 
few survivors at extreme ages are an intrinsically healthier subset of 
the initial population) and temporal heterogeneity (the aging 
process itself slows down in time).

One might have expected, in a field as thoroughly ploughed over as 
the asymptotic behavior of Markov processes, that any problem thrown 
up by applied science would already have been effectively answered in 
the probability literature.  This problem seems to be an exception.  In 
particular, one-dimensional diffusions which are killed during their 
motion appear to have been overlooked.  The case of a diffusion with one 
inaccessible boundary and one regular boundary, at which it is killed (or 
partly killed and partly reflected), with no internal killing, has 
been elucidated by P. Mandl \cite{pM61}.  His methods were further developed, in the more specialized case of pure killing at the boundary, by P. Collet, S. Mart\'inez, P. Picco and J. San Mart\'in \cite{CMSM95,MPSM98,MSM01,MSM04}.

Some of the general principles are the same as those for countable discrete-time Markov chains, first discussed in \cite{aY47} (in the context of branching processes), and further developed by D. Vere-Jones alone \cite{dVJ62} (where the R-theory, later developed by R. Tweedie \cite{rT74}, was first introduced) and with E. Seneta \cite{SV66}.    The corresponding theory for continuous-time countable-state processes was introduced by J. F. C. Kingman \cite{jK63}.  We note that the results of \cite{CMSM95} are presaged in the work of E. van Doorn \cite{eVD91} and of M. Kijima and E. Seneta  \cite{KS91} on birth-death chains conditioned on not having dropped below 0, which itself is built upon a birth-death-chain countable version of Mandl's work by J. Cavender \cite{jC78}.  The study of quasistationary distributions for birth-death chains was initiated by P. Good \cite{pG68}.  There is a strong analogy between birth-death chains killed at 0 and diffusions killed at 0.  Some features of the problem, in particular the role of the initial distribution, are easier to understand in the birth-death chains, and are well described by \cite{eVD91} in particular.  A slight generalization is offered by \cite{FKMP95}, where the Markov process in the positive integers is allowed to have arbitrary jumps, but is still constrained to have killing only at 0, and to take a long time to reach 0 from far away.  Similarly, quasistationary distributions for a particular density-dependent branching process were studied by G. H\"ogn\"as \cite{gH97}.

More general killing has primarily been considered in the finite state-space setting.   General theory for quasistationary distributions on finite state-spaces is worked out in \cite{DS65}, and at greater length in \cite{eS73}; applications to genetics are discussed in \cite{eS66}.  In this context, it is not hard to show that quasistationary distributions exist (under basic irreducibility assumptions), and also that the distribution of the process conditioned on survival converges to that distribution (sometimes referred to as the ``Yaglom limit''). Only recently have comparable results been developed for Markov chains on infinite state spaces, by F. Gosselin \cite{fG01} (for countable state spaces) and by J. B. Lasserre and C. E. M. Pearce \cite{LP01} (for general state spaces).   Both of these have results similar to those of the present work (allowing for the substantial technical difference in setting), but also fundamentally different, in that they rely on a Lyapunov-type criterion.  We are not aware of any earlier work that derives these kinds of convergence results for diffusions with internal killing.

In the present work, we extend the results of \cite{pM61} and \cite{CMSM95} to the case of general killing, and prove new results which would not arise when killing occurs only at 0.  The essential result is a dichotomy: The process either escapes to $\infty$ with probability 1 (so that the probability of finding it in any compact interval goes to 0), or it settles into the quasistationary distribution given by the top eigenfunction of the adjoint generator.  Along the way, we shore up the foundations of this subject, supplying details which were left vague.\footnote{As we discuss in the context of the proofs, some details of P. Mandl's paper needed to be fleshed out, to show convincingly that the methods could be generalized.  In addition, some fundamental results about the generators of diffusions with killing had not been written down.}

\renewcommand{\thefootnote}{+}

Our methods are more ``hands-on'' than those of
\cite{CMSM95}.  In part, this is compelled by
a serious gap in their proof.  In their Lemma 2 it is stated that a
sequence $t_{n}\to\infty$ may be found such that the measures $\P_{h}
\ls X_{t_{n}} \in \cdot \cond \tp>t_{n}\rs$ converge, where $h$ is an
initial distribution with compactly supported density.  Their proof
simply attributes this fact to compactness of the space of probability
measures on the extended real half-line.  While this is, of course,
correct for any single $h$, they claim, with no further justification,
that the same limit measure may be obtained for an arbitrary countable
set of initial distributions.  This is equivalent, in the notation
that we introduce in Lemma \ref{L:omegat}, to claiming that
$F(\nu,\R^{+})$ is independent of $\nu$.  As we show in Proposition
\ref{P:collet}, this is indeed equivalent to the desired result.
While we cannot rule out the possibility that this fact may be proved
in some fashion trivial enough to be neglected, we have been unable to
find a direct proof.\footnote{The authors of \cite{CMSM95} have provided an alternative proof of their result, in a private communication, for their special case where killing occurs only at 0.  Their proof of the lemma amounts essentially to a direct proof of their theorem, by proving Mandl's extra condition \eqref{E:extracond1}.  Indirectly, this also proves their Lemma 2.  There is, as yet, no direct proof of their Lemma 2, hence no reason to think that the method proposed in \cite{CMSM95} can be made to work, even in the case of no internal killing.}

Where we have made good use of a generalized version of the constructions from \cite{CMSM95}, is in deriving our Theorem \ref{T:limits}, a result that has no parallel in their setting, but can only arise in
the case of nonzero internal killing.  If the killing rate is faster at infinity than near the origin, then escape to infinity is
excluded: any mass that moves away from the origin will become
relatively minuscule over time.  This allows us to determine, when the limit of the killing rate at infinity is greater than the quasistationary killing rate, that the process does not run off to infinity.

\section{Definitions and assumptions}  \label{sec:defout}
\nc{\ulm}{\underline{\gl}}
\subsection{The process} \label{sec:proc}
We consider a one-dimensional diffusion on the interval 
$(0,r)$, where $r$ is finite or infinite.  We will always assume 
that 0 is a regular boundary.  If the right boundary is $\infty$, we 
assume it to be inaccessible (that is, natural or entrance); if it is 
finite, we assume it to be regular.  When $r$ is not explicitly mentioned, it is assumed to be $\infty$.

The diffusion rate $\gs$ is 
assumed continuously twice differentiable.  The drift $b\in C^{1}((0,r))$ is
assumed continuously differentiable.  The killing rate, denoted 
$\kappa$, is only required to be continuous and nonnegative.   We write $\tau_{x}$ for the first time the process hits state $x$.  The time at which the process is killed, transitioning to the cemetery state $\partial$, is thus represented by $\tp$.

   We may simplify the problem somewhat by assuming that $\gs$ is           
   identically 1.  There is no loss of generality since (following          
   \cite{wF52}) we can replace $X_{t}$ by $Y_{t}=F(X_{t})$, where           
   $F(x):=\int_{x_{0}}^{x}du/\gs(u)$, where $x_{0}$ is an arbitrary point 
   in $(0,r)$.  The killing rate
for $Y_t$ becomes $\kappa(F^{-1}(Y_t))$, while the drift may be
computed by It\^o's formula to be
   \begin{equation}
     \label{E:newdrift}
     \frac{b(F^{-1}(Y_{t}))}{\gs(F^{-1}(Y_{t}))}-\gs'(F^{-1}(Y_{t})).
   \end{equation}
   Since $\gs$ is positive on the interior of the interval, this 
   transformation is finite.  From now on, we will always assume, 
   unless otherwise stated, that the diffusion rates of our 
   one-dimensional diffusions are identically 1.  We note that in this normalization the speed measure has density $e^{B(x)}$ with respect to Lebesgue measure, where $B(x):=2\int_{0}^{x}b(z)dz$, while the scale function is $s(z)=\int_{x_{0}}^{z} e^{-B(x)}dx$.  We denote by $\beta$ the measure on $\R^{+}$ that has density $e^{-B(x)}$ with respect to Lebesgue measure.

The diffusion semigroup acts on the Banach space $C([0,r])$ of functions continuous on the closed interval $[0,r]$, with the supremum norm.  (When $r=\infty$ this is the extended half-line.  We will use the notations $[0,\infty)$ and $\R^{+}$ interchangeably for the half-line including 0.)  As a strongly continuous contraction semigroup, it has an infinitesimal generator, which is the closure of the operator
\begin{equation}
    \mL \phi =\oh \phi''(x)  +b(x) \phi'(x)-\kappa(x)\phi(x),
    \label{eq:opnokill}
\end{equation}
acting on a subspace of the twice differentiable functions whose image under the operator is also in $C([0,r])$.  (The subspace is determined by lateral conditions at the endpoints.  Each regular endpoint requires a boundary condition to restrict the generator; the boundary conditions restrict the domain of the pregenerator, which is then closed to make a generator.)  The dual space is a subspace of $\eu{BV}([0,r])$, the signed measures with bounded variation on the closed interval $[0,r]$, with total-variation norm.  There is a natural inclusion of $L^{1}$ into $\eu{BV}((0,r))$, taking $f\in L^{1}$ to the signed measure with density $f$ on $(0,r)$.  This identifies $L^{1}$ with the subspace of absolutely continuous signed measures.  Since the range of the resolvent is contained in this subspace, the semigroup may be restricted to $L^{1}$.  (This is proved in Theorem 13.3 of \cite{wF52}, which we extend to diffusions with internal killing in Lemma \ref{L:gen}.  Because the space $\eu{BV}([0,r])$ includes masses at the endpoints, in general these could be retained in the image.  In fact, though, when the boundaries are natural, entrance, or regular without ``stickiness'', point masses at the endpoints are excluded from the image of the resolvent.)  The adjoint generator is the closure of the operator
\begin{equation}
    \mL^{*} \phi =\oh \phi''(x) -\lp b(x) \phi(x) \rp'-\kappa(x)\phi(x),
    \label{eq:opnokilladj}
\end{equation}
acting on the space of twice differentiable functions whose image is in $L^{1}$, with boundary conditions implied by the nature of the endpoints.  (This description of the generator and the adjoint is given by W. Feller in \cite{wF52}.)

The generator describes the behavior of the process in the interior of
the state space.  The behavior at the endpoints is determined by a
boundary condition on the domain.  If 0 is a regular boundary, we need a boundary condition to describe the behavior.  We assume that there is no ``holding'' or ``stickiness'' at the boundary --- that is, the speed measure $m$ satisfies $m(\{0\})=0$ --- and we assume that $b$ is continuous at a regular boundary.  For this case, Theorem 19.3 of \cite{wF52} gives the boundary condition on the domain of $\mL^{*}$ as
\begin{equation}  \label{E:fellerbound}
(1-p_{0}) \phi(x) = p_{0}\lp \oh\phi'(0)-b(0) \phi(0)\rp,
\end{equation}
where $0\le p_{0}\le 1$ is an arbitrary parameter.  (Note: Everywhere that we have $\oh\phi'-b\phi$ for the boundary conditions, the sources we cite have twice that.  Of course, this is arbitrary, since the normalization of $p_{0}$ is arbitrary; but it seems appropriate, because the sources are normalized with $\gs^{2}=2$, whereas we have $\gs^{2}=1$.)
It is a well-known fact that $(1-p_{0})/p_{0}$ represents the rate of absorption at 0.  
If $r$ is a regular boundary as well, the analogous condition holds there:
\begin{equation}  \label{E:fellerboundr}
(1-p_{r}) \lim_{x\to r}e^{-B(x)}\phi(x) = -p_{r}\lim_{x\to r} \lp \oh\phi'(x)-b(x) \phi(x)\rp.
\end{equation}
If $r$ is an entrance boundary, there is the simpler condition
\begin{equation}  \label{E:fellerboundentrance}
\lim_{x\to r}\lp \oh\phi'(x)-b(x) \phi(x)\rp=0.
\end{equation}
(These boundary conditions are given in \cite{wF52} for the forward semigroup and the adjoint, in the absence of internal killing. In \cite{IM65} the same conditions are derived for general killing, but only for the forward semigroup.  We transfer this result to the adjoint in Lemma \ref{L:gen}.)

The adjoint semigroup is not discussed in \cite{IM65}.  The arguments
of \cite{wF52} carry over directly.  We recapitulate them here,
because an exact description of the adjoint generator --- and of its
domain, in particular --- will be indispensable to the proof of
Theorem \ref{T:mainnatreg}.  We define $\mC$ to be
the set of twice continuously differentiable functions on $\R^{+}$
satisfying the boundary condition \eqref{E:fellerbound} at 0, and the
appropriate analogue condition at $r$ in the case of two regular
endpoints.  Given an operator $A$ and a strongly continuous
contraction semigroup $T_{t}$, we will write ``$A$ generates $T_{t}$''
or ``$T_{t}$ is generated by $A$'' as shorthand for ``$A$ is closable,
and the closure of $A$ is the generator of $T_{t}$''.

There are several additional conditions that we will need to impose on
our diffusions.  To begin with, we need
\begin{equation*}
\tag{LP}  
\mL^* \text{ is in the limit-point case at }\infty.
\end{equation*}
(For an account of the theory of singular Sturm-Liouville problems,
and the concepts of ``limit-point'' and ``limit-circle'' in
particular, see either chapter 9 of \cite{CL55} or chapter 5 of
\cite{kY60}.)  Our main results will all assume condition (LP).  This is fairly
technical, but there is a simple condition on the diffusion parameters
which suffices.

\begin{Lem}
\label{L:LP}
The condition
\begin{equation*}
        \liminf_{z\to\infty} z^{-2} \lp b(z)^{2} + b'(z) + 2\kappa(z) 
        \rp >-\infty\quad \text{or}\quad \kappa \equiv 0.\tag{LP'}
    \end{equation*}
implies (LP).
\end{Lem}

\begin{proof}
  When $\kappa\equiv 0$, the proof is given in \cite{pM61}.  The limit-point property follows directly from the fact that
  $\phi(z)=e^{B(z)}\int_{0}^{z}e^{-B(x)}dx$ is a solution to $\mL^{*}
  \phi=0$. The assumption that $\infty$ is inaccessible implies that
  this is not in $L^{2}(\beta)$.
 
We make the standard transformation
\begin{align*}
\hat{\mL}\phi(x) &:= e^{-B(x)/2} \mL \lp e^{B(x)/2}\phi(x)\rp\\
&=\oh \phi''(x) -\oh\lp b(x)^{2}+b'(x)+2\kappa(x)\rp \phi(x)
\end{align*}
The eigenfunctions of $\hat{\mL}$ may be transformed into those of 
$\mL^*$ simply by multiplying by $e^{B(x)/2}$.  Thus $\hat{\mL}$, which 
is self-adjoint with respect to Lebesgue measure, is in the 
limit-point case at $\infty$ if and only if this is true of $\mL^*$, which is 
self-adjoint with respect to $e^{B(x)}dx$. But Corollary 9.2.2 of 
\cite{CL55}, together with the assumption (LP'), 
implies that $\mL^{*}$ is in the limit-point case at $\infty$.
\end{proof}

Note that (LP') is not an onerous condition: it can fail only if $-b'(z)$ becomes arbitrarily large (as $z\to\infty$), while $b(z)^{2}$ remains relatively small.

\subsection{Eigenfunctions of the generators}
For any $\gl$, we define $\phi_{\gl}$ to be the unique solution to the
initial value problem
\begin{equation}  \label{E:ivp}
\mL^* \phi_\gl =-\gl \phi_{\gl}\quad
 \text{with }  \phi_{\gl}(0)=p_{0} \text{ and } 
\oh \phi'_{\gl}(0)- b(0)\phi_{\gl}(0)=1-p_{0}.
\end{equation} 
We define 
\begin{equation}  \label{E:ulm}
\ulm=\sup\ls \gl \, :\, \phi_{\gl}(x) \text{ is nonnegative for all }x \rs.
\end{equation}
Similarly, we define $\psi_{\gl}$ to be the unique nonzero solution to
the initial value problem
\begin{equation}  \label{E:ivp2}
\mL \psi_{\gl} =-\gl \psi_{\gl}\quad
 \text{with }  \psi_{\gl}(0)=p_{0} \text{ and } \oh\psi'_{\gl}(0)=1-p_{0}.
\end{equation}
Note that $\psi_{\gl}(x)=e^{-B(x)}\phi_{\gl}(x)$.  If $\psi$ is any
$\gl$-invariant function --- that is, a function $\psi$ such that
$P_{t}\psi=e^{-\gl t}\psi$ --- then $\psi$ is a multiple of
$\psi_{\gl}$.  The converse need not be true: that is, $\psi_{\gl}$
need not be $\gl$-invariant.  We follow Mandl in defining
$\gl_{\varrho}$ to be the smallest point of increase of the spectral
measure $\varrho$ for the operator $\mL^{*}$.  (We also follow Mandl
in reversing the sign on the spectral measure; where the spectrum
should be concentrated on $(-\infty,-\gl_{\varrho}]$ --- which would
correspond to considering eigenvalues $\phi_{\gl}$ with eigenvalue
$\gl$, where $\gl$ is negative --- we have let $\phi_{\gl}$ have
eigenvalue $-\gl$, and put $\varrho$ on $[\gl_{\varrho},\infty)$.)

We follow Mandl in defining $\Psi$ to be
the unitary map from $L^{2}(\beta)$ to $L^{2}(\varrho)$ defined by the
eigenfunctions.  That is, for $f\in L^{2}(\beta)$,
\begin{equation}  \label{E:definePsi}
\Psi f (\gl) = \lim_{M\to\infty} \int_{0}^{M} f(x)\phi_{\gl}(x) e^{-B(x)}dx,
\end{equation}
where the limit is understood in the $L^{2}(\beta)$ sense.  The
existence of this limit, and the fact that it is unitary, are shown in
Theorem 3.1 of \cite{CL55}.  The measure $\beta$ does not appear
there, but the essential assumption is that the differential operator
is self-adjoint.  Since our differential operator $\mL^{*}$ (acting on
the appropriate domain $\hat{\eu{D}}$, defined in section
\ref{sec:compact}) is self-adjoint with respect to $\beta$, the result
carries over.

An essential step in Mandl's method is his Lemma 2, the fact that
$\ulm=\gl_{\varrho}$.  This depends only on the assumption that the
operator is in the limit-point case at $\infty$.  The proof is
sketchy, though, in the original, which might leave one in some doubt
of whether it can properly be extended to our more general setting.

\begin{Lem}
\label{L:bottom}
Suppose (LP) holds.  Then $\ulm=\gl_{\varrho}$ is finite, and the set
of $\gl$ such that $\phi_\gl$ does not change sign is
$(-\infty,\ulm]$.  Similarly, the set of $\gl$ such that $\psi_\gl$ does not change sign is $(-\infty,\ulm]$.
\end{Lem}

\begin{proof}
Let $\mI$ be the set of $\gl$ such that the boundary
value problem \eqref{E:ivp} on $[0,\infty)$ has a nonnegative (and nontrivial)
solution.  Theorem 8.2.1 of \cite{CL55} tells us that the
$r$-truncated boundary-value problem has solutions for infinitely many
$\gl$.  By the $r$-truncated problem we mean the initial-value problem
on $[0,r]$ with the same condition at 0, and the additional condition
$\phi(r)=0$.  These $\gl$ may be ordered as $0\ge
-\gl_0>-\gl_1>\cdots$, where $\lnf \gl_n=\infty$, and the
eigenfunction $\phi_{\gl_n}$ has exactly $n$ zeroes on $(0,b)$.  Any
$\gl$ for which there is a solution to the $r$-truncated problem for
some $r$ is automatically not in $\mI$.  This implies that the
complement of $\mI$ is unbounded above.

\nc{\ulam}{\gl_+} 
\nc{\olam}{\gl_-} 

Suppose now that we had three values $\ulam<\gl^*<\olam$, with
$\ulam,\olam\in\mI$, but $\gl^*\notin\mI$.  Since the zeroes
in any solution to the differential equation are discrete, we may
define an extended real-valued function $\zeta(\gl)$ on
$[\ulam,\olam]$ as the smallest zero of $\phi_\gl$, or $\infty$ if the
function has no zero.  On any compact set of $z$ and $\gl$,
$\phi_\gl(z)$ is uniformly continuous in $z$ and $\gl$ (by Theorem
1.7.5 of \cite{CL55}).  This implies that $\zeta$ is continuous. 
Since $\zeta(\ulam)=\zeta(\olam)=+\infty$, and since $\zeta(\gl^*)$ is
finite, there must be $\gl'<\gl''$, both in $[\ulam,\olam]$, such that
$\zeta(\gl')=\zeta(\gl'')$.  But this would imply that the
$b$-truncated problem for $b=\zeta(\gl')$ has a nonnegative solution
for two values of $\gl$, namely $\gl'$ and $\gl''$.  This contradicts
the above-mentioned Theorem 8.2.1 of \cite{CL55}.  The same argument
shows that we cannot have $\ulam,\olam\notin\mI$, but
$\gl^*\in\mI$.

We may conclude that $\mI$ is either empty or an interval unbounded below.  
Mandl shows clearly that $\gl_\varrho\in\mI$.  The proof depends only on
the fact that the initial-value problem is in the limit-point case at
$\infty$, which we have assumed.  On the other hand,
$\gl_\varrho$ is a point of increase for $\varrho$ defined as
a limit as $r\to\infty$ of step functions with jumps at the
eigenvalues for the $r$-truncated problem.  There must be
$\gl>\gl_\varrho$ arbitrarily close, such that $\gl$ is an
eigenvalue for the $r$-truncated problem for some $r$.  These $\gl$
are not in $\mI$.  Thus, $\mI=(-\infty,\gl_\varrho]=(-\infty,\ulm]$.

The same argument works for the boundary-value problem \eqref{E:ivp2}.
\end{proof}

\subsection{Asymptotic properties}
We say that the process $X_{t}$ has an {\em asymptotic killing rate $\eta$ when started from the probability distribution $\nu$} if
\begin{equation}  \label{E:sakr}
\ltf \P_{\nu} \ls \tp>t+s \cond \tp>t \rs = e^{-\eta s}.
\end{equation}
for any positive $s$.  A trivial consequence is that
\begin{equation}  \label{E:akr}
\eta=-\ltf t^{-1}\log \P_{\nu}\ls \tp >t\rs.
\end{equation}

If $\phi\in C_{0}([0,\infty))$ is nonnegative, we say that the process $X_{t}$ {\em converges from the initial distribution $\nu$ to the quasistationary distribution $\phi$ on compacta} if for any positive $z$, and any Borel $A\subset [0,z]$,
\begin{equation}
    \ltf \P_{\nu}\ls X_{t}\in A \cond X_{t}\le z \rs = 
    \frac{\int_{A}\phi(z)dz}{\int_{0}^{z}\phi(z)dz}.
    \label{E:firstlim2}
\end{equation}
It should be noted that this property, while interesting, is far short of what one would like to prove.  In particular, it leaves open the possibility that most of the mass --- an amount which falls off at a rate strictly slower than $e^{-\ulm t}$ --- runs off to $\infty$, even while the part of the distribution which stays below $z$, for any fixed $z$, converges to the distribution with density $\phi_{\ulm}$.
We say that $X_{t}$ {\em converges from the initial distribution $\nu$ to the quasistationary distribution $\phi$} if \eqref{E:firstlim2} holds with $z=\infty$; that is, for any Borel subset $A\subset [0,\infty)$,
\begin{equation}
    \ltf \P_{\nu} \ls X_{t}\in A \cond \tp>t \rs = 
    \frac{\int_{A}\phi(z)dz}{\int_{0}^{\infty}\phi(z)dz}.
    \label{E:firstlim3}
\end{equation}

We say that $X_{t}$ {\em escapes to infinity} from $\nu$ if
\begin{equation}  \label{E:escape}
\lim_{t\to\infty} \P_{\nu} \ls X_{t}\le z \cond \tp>t \rs =0
\end{equation}
for all $z\in \R^{+}$.

An elementary consequence of the definitions is:

\begin{Prop}  \label{P:esc}
If $X_{t}$ started at $\nu$ converges to the quasistationary distribution $\phi$ on
  compacta, and $\izf \phi(y)dy=\infty$ then $X_{t}$ started at $\nu$ escapes to
  $\infty$.
\end{Prop}

\begin{proof}
Given any $z'>z>0$,
\begin{align*}
\lstf \P_{\nu}\ls X_{t}\le z \cond \tp>t\rs &\le 
   \lstf \P_{\nu}\ls X_{t}\le z \cond X_{t}\le z'\rs\\
   &= \frac{\int_{0}^{z} \phi(y)dy}{\int_{0}^{z'} \phi(y)dy}.
\end{align*}
Sending $z'$ to $\infty$ proves the result.
\end{proof}

\nc{\lzn}{\limsup_{z\to\infty}\limsup_{n\to\infty}}
\nc{\llm}{\tilde{\gl}}

\subsection{The initial distribution}  \label{sec:id}
The convergence properties may depend on the initial distribution.  In
particular, if the tail of the initial distribution is too heavy,
there may be mass wandering gradually in from the tails more slowly
than the decay rate of the mass starting on any compact set.  (In the case of countable-state discrete Markov chains, conditions on the initial state leading to convergence to the quasistationary distribution are discussed in \cite{SV66}.) 
We adopt the following sufficient condition from Mandl:
\begin{equation*}
\text{The initial distribution has a density } f\in L^{2}(\beta) \text{ with }
\liminf_{\gl \downarrow \ulm} \Psi f(\gl) >-\infty.
\tag{ID}
\end{equation*}
(Unless otherwise indicated, densities are always meant to be taken with respect to Lebesgue measure.)
Since we are concerned with asymptotic properties, there is no change
if we replace $X_{t}$ by $X_{t+s}$, for some positive $s$.  Thus, we
may substitute
\begin{equation*}
\begin{split}
\text{If }X_{0}\text{ has distribution }&\nu,\text{ then }\te s\ge 0 \text{ for which the distribution of }X_{s} \\
&\text{ has density } f\in L^{2}(\beta), \text{ with }\liminf_{\gl \downarrow \ulm} \Psi f(\gl) >-\infty.
  \end{split}
\tag{ID'}
\end{equation*}

We will use (IDC) to denote the condition that a distribution satisfies (ID') and has compact support.  The following is a trivial consequence of the definition \eqref{E:definePsi} of $\Psi$:
 
\begin{Prop}  \label{P:ID}
If the initial distribution has a density with compact support, then condition (ID) is satisfied.
\end{Prop}

It seems plausible that condition (ID') should be satisfied for any initial distribution with compact support (which seems to be implied by comments in \cite{pM61}), but we have not been able to show this.

\subsection{The $\omega$ function}  \label{sec:omega}
We define the spaces $\mc$, comprising probability distributions on $[0,\infty)$ with compact support, and $\mo$, comprising subprobability distributions on $\R^{+}$.  When $p_{0}=0$ (that is, in the case of complete killing at 0), we exclude the point $\delta_{0}$ from $\mc$.  We supply $\mo$ with the topology of weak convergence on compact sets, while $\mc$ has the stronger topology in which $\nu_{n}\to\nu$ if and only if the sequence converges weakly, and $\bigcup_{n}\supp\nu_{n}$ is bounded.  (That is, the open sets are generated by weak open sets, together with $\ls \nu\, : \, \supp \,\nu \subset (0,R)\rs$ for $R>1$.)

We imitate \cite{CMSM95} in defining functions
\begin{equation}  \label{E:omega}
\omega_{t}(\nu):= \frac{\P_{\nu}\{\tp>t \}}{\P_{1}\{\tp>t \}} \text{ for } \nu\in\mc,
\end{equation}
mapping compactly supported probability distributions on $\R^{+}$ to $\R^{+}$.  We also define
\begin{equation}  \label{E:ft}
F_{t}(\nu,\cdot):= \P_{\nu}\ls X_{t}\in \cdot \cond \tp > t \rs .
\end{equation}
When $\nu=\delta_{x}$, we will also write $F_{t}(x,A)$ and $\omega_{t}(x)$.  Note that 
\begin{equation}  \label{E:omeganu}
\omega_{t}(\nu)=\izf \omega_{t}(x) d\nu(x)
\end{equation}
for any compactly supported $\nu$.  Since $\omega_{t}(\nu)$ is nonzero (except in the case of complete killing at 0, when $\nu=\delta_{0}$), we also have, for all $\nu\ne \delta_{0}$,
\begin{equation}  \label{E:Fnu}
F_{t}(\nu,\cdot)=\frac{\izf F_{t}(x,\cdot) \omega_{t}(x) d\nu(x)}{\izf \omega_{t}(x) d\nu(x)}.
\end{equation}
We note that a version of $\omega_{t}$ (where the argument was confined to point masses) was introduced in \cite{CMSM95}.

Finally, we define
\begin{equation}  \label{E:definea}
a_{t}(\nu,r):=\P_{\nu} \ls \tp>t+r\cond \tp>t \rs = \int F_{t}(\nu,dy) \P_{y} \ls \tp>r \rs.
\end{equation}

For some results, we will need the following condition:
\begin{equation*}  \tag{GB}
\fa s\ge 0, \quad \E_{x} \sup_{t\ge 0}\omega_{t}(X_{s})<\infty
\end{equation*}
This seems to be an abstract condition, which might be hard to verify.  In fact, though, it is a simple consequence of some fairly weak conditions on the diffusion parameters $b$ and $\kappa$.

Consider first, the situation that $\kappa$ is monotonically increasing and there is complete reflection at 0 (as in our example of section \ref{sec:example}).  We can couple a copy of the process $X_{t}$ started at 1, to a copy $X'_{t}$ started at $x\ge 1$ in such
a way that $X_{t}\le X'_{t}$ while both processes are alive.  At any time when both processes are alive, 
the killing rate must be larger for $X'$ than for $X$, 
so we can also arrange the coupling so that $X'$ is killed before $X$.  This implies that $\omega_{t}(x)\le 1$ for all $x$ and all $t$,  
and hence condition (GB) is trivially satisfied.

If $\kappa$ is not monotone, we need to restrict its growth, as well as the growth of the negative drift.
We note that there is a general bound on $\omega_{t}(x)/\omega_{t}(y)$, obtained by factoring $\P_{y}\{\tp>t \}$ over the time (if any) when the process started at $y$ first hits $x$.  We have
\begin{align*}
\frac{\omega_{t}(x)}{\omega_{t}(y)}=\frac{\P_{x}\{\tp > t\}}{\P_{y}\{\tp > t\}}
&\le\frac{\P_{x}\{\tp > t\}}{\int_{0}^{t} \P_{x}\{\tp > t-u\}\P_{y}\{\tau_{x}\in du\}}\\
&\le \P_{y} \ls \tau_{x}<\tp \rs^{-1}=:\omega_{*}(x,y)
\end{align*}
This yields
\begin{align}  \label{E:unifomegabd}
\omega_{*}(y,x)^{-1} \le\frac{\omega_{t}(x)}{\omega_{t}(y)}\le&\omega_{*}(x,y),\text{ and}\\
\left| \frac{\omega_{t}(x)}{\omega_{t}(y)}-1\right|\le \omega_{**}(x,y)&:= \omega_{*}(x,y)-\omega_{*}(y,x)^{-1}.\label{E:unifomegabd2}
\end{align}
For \eqref{E:unifomegabd2} we use the fact that if $A\le B\le C$ and $A\le 1\le C$, then $|B-1|\le C-A$.

\begin{Lem}
\label{L:omegabounded}
Suppose there are positive constants $\kappa_{*}$, $b_{*}$, $b_{**}$, and $\beta$ such that either
\begin{equation*}  \tag{GB'}
\lv b(y)\rv\le b_{*}y \quad \text{and}\quad \kappa(y)\le \kappa_{*}y
\text{ for }y \text{ sufficiently large}
\end{equation*}
or
\begin{equation*}  \tag{GB''}
-b_{*}y^{\beta}\le b(y)\le b_{*}y,\quad b'(y)\le b_{**}y^{2}\quad\text{and}\quad \kappa(y)\le \kappa_{*}y
\text{ for }y \text{ sufficiently large}.
\end{equation*}
Then for all $s,t,x\in\R^{+}$,
\begin{equation}  \label{E:omegabounded}
 \E_{x} \sup_{t\ge 0}\omega_{t}(X_{s}) \le \E_x\omega_{*} (X_{s},1) <\infty.
\end{equation}
\end{Lem}

\nc{\tl}{\tp^{(\gl)}}

\begin{proof}
Integrating by parts, we need to show that
$$
\lim_{w\to\infty} \frac{\P_{x}\{X_{s}\ge w \}}{\P_{1}\{\tau_{w}<\tp\}} \, - \,
\int_{1}^{\infty} \frac{\P_{x}\{X_{s}\ge z \}}{\P_{1}\{\tau_{z}<\tp\}^{2}} \left(
\frac{d\phantom{z}}{dz} \P_{1}\{\tau_{z}<\tp\} \right)  dz 
$$
is finite.
For $\gl$ nonnegative, define $\tl$ to be the minimum of $\tp$ and an independent exponential stopping time with expectation $1/\gl$.  Then
$$
\P_{x}\ls X_{s}\ge z \rs\le \P_{x} \ls \tau_{z}< s \rs \le e^{\gl s} \P_{x} \ls \tau_{z}<\tl \rs.
$$
Since $\tl\le \tp$, it suffices to show
\begin{equation}  \label{E:sufftoshow}
- \int_{x}^{\infty} \frac{\P_{x}\{\tau_{z} <\tl \}}{\P_{1}\{\tau_{z}<\tp\}^{2}} 
\left( \frac{d\phantom{z}}{dz} \P_{1}\{\tau_{z}<\tp\} \right)dz<\infty.
\end{equation}

Observe now that $\P_{y} \ls \tau_{z}<\tl \rs$ for $y\in(0,z)$ is the Laplace transform of $\tau_{z}$.  On this interval it coincides with eigenfunction $\psi_{-\gl}(y)/\psi_{-\gl}(z)$ of $\mL$.  What we need to show, then, is that
\begin{equation}  \label{E:reallyneed}
\int^{\infty} \frac{-\psi'_{0}(z)}{\psi_{-\gl}(z)} dz <\infty
\end{equation}
for some positive value of $\gl$.

We define
\begin{equation}  \label{E:smp}
g_{\gl}(x):=\frac{\psi'_{-\gl}(x)}{\psi_{-\gl}(x)}.
\end{equation}
Note that $g_{\gl}$ satisfies
\begin{equation} \label{E:ggl}
g_{\gl}(x)=\lim_{\Delta x\downarrow 0}
(\Delta x)^{-1}\P_{x-\Delta x}\ls \tau_{x}>\tl \rs,
\end{equation}
from which it follows that $g_{\gl}(x)$ is positive, and nondecreasing in $\gl$.
With this notation, the condition \eqref{E:reallyneed} becomes
\begin{equation}  \label{E:reallyreallyneed}
\int^{\infty} g_{0}(z) \exp\left\{ \int_{1}^{z}\lp g_{0}(y)-g_{\gl}(y)\rp \right\} dz <\infty.
\end{equation}

\nc{\ok}{\overline{k}}
\nc{\ob}{\overline{b}}
We now need an upper bound on $g_{\gl}(x)$.  Fix any positive $\gd$.  We start the diffusion at $x-\gep$, where $\gep$ is assumed smaller than $\gd$. Note that if we increase the killing rate on $(x-\gd,x)$ to $\ok:=\sup_{y\in (x-\gd,x)}\kappa(y)$, and kill the process when it reaches $x-\gd$, that increases $g_{\gl}(x)$, since the probability of being killed before reaching $x$ is increased.  If we then decrease the drift rate on $(x-\gd,x)$ to the constant $\ob:=\inf_{y\in (x-\gd,x)}\kappa(y)$, that further increases $g_{\gl}(x)$.  We can then apply equation 2.3.0.5 of \cite{BS96} to compute
$$
g_{\gl}(x)\le \frac{d\phantom{x}}{dx}\left(-e^{\ob x}
\frac{\sinh(\ga(\gd-x))}{\sinh(\ga \gd)}\right)_{x=0}
\le -\ob+\ga \coth(\ga \gd) \le -2(\ob\wedge 0)+\sqrt{2\ok+2\gl}+\gd^{-1}
$$
where $\ga=\sqrt{\ob^{2}+2\ok^{2}+2\gl}$.  If (GB') holds, we fix, arbitrarily, $\gd=1$, and we see that for $x$ sufficiently large
\begin{equation}  \label{E:glbound0} 
g_{\gl}(x)\le 2b_{*}x+\sqrt{2\kappa_{*}}x+\sqrt{2\gl}+1.
\end{equation}  
If (GB'') holds, we take $\gd=1/x$.  Then for $x$ sufficiently large
\begin{equation}  \label{E:glbound1}
g_{\gl}(x)\le -2(b(x)\wedge 0)+2b_{**}x+\sqrt{2\kappa_{*}}x+\sqrt{2\gl}+x.
\end{equation}

Define $h(y):=g_{\gl}(y)- g_{0}(y)$, and $k(y):=g_{\gl}(y)+g_{0}(y)+2b(y).$  
Using either \eqref{E:glbound0} or \eqref{E:glbound1}, we see that
\begin{equation}  \label{E:glbound2}
k(y)\le Cy+\sqrt{2\gl}+C',
\end{equation}
where $C$ and $C'$ are constants independent of $\gl$.

From the equation \eqref{E:ivp2} we derive the relation
\begin{equation}
h'(y)=2\gl - k(y) h(y) .\label{E:heq}
\end{equation}
This equation has the general solution
\begin{equation}  \label{E:hsol}
h(z)=h(1)e^{-\int_{1}^{z}k(y)dy} +2\gl \int_{1}^{z} e^{-\int_{y}^{z}k(x)dx}dy.
\end{equation}
Let $C''=\max\{\sqrt{2\gl}+C',1\}$.  By \eqref{E:glbound2}, for $z\ge 1$,
\begin{align*}
h(z)&\ge 2\gl \int_{z-1/C''z}^{z} \exp\Bigl\{-C(z^{2}-y^{2})/2-C''(z-y)\Bigr\}dy\\
&\ge (\sqrt{2\gl}-C') e^{-2C-1}z^{-1}.
\end{align*}

Thus,
$$
\exp\left\{ \int_{1}^{z}\lp g_{0}(y)-g_{\gl}(y)\rp \right\}\le z^{-(\sqrt{2\gl}-C') e^{-2C-1}}.
$$
Since $g_{0}$ grows at a polynomial rate, this proves \eqref{E:reallyreallyneed} for $\gl$ large enough.
\end{proof}

\section{Outline of the paper and statement of main results}  \label{sec:main}
We will consider here mainly the case where 0 is a regular endpoint,
and $r=\infty$ is an inaccessible endpoint.  The case of two regular
endpoints is an application of basic results about compact operators,
and differs only slightly from the standard $\kappa\equiv 0$ case,
such as is presented in \cite{rP85}.  We present the more general
proof only briefly. The case of two inaccessible endpoints also would
offer no additional novelty, except to complicate the notation.

Mandl's argument for diffusions on $\R^{+}$ proceeds in two steps.  In
the first step, he proves that $X_{t}$ converges to the
quasistationary distribution $\phi_{\ulm}$ on compacta.  In order to
remove the restriction to compact sets, Mandl imposes a condition
equivalent to
\begin{equation}
        \label{E:extracond1}
    \int_{0}^{\infty} e^{B(z)} dz <\infty \text{ and }
    \lim_{z\to\infty}\int_{\underline{\lambda}}^\infty
     \left( \int_z^\infty \phi_\lambda (y) dy \right)^2 \varrho(d\lambda) = 0,
    \end{equation}
    implying a strong drift away from $\infty$.  Collet {\em et al.} \cite{CMSM95}
    offer a different approach to removing the restriction to finite
    $z$.  Under the assumption that the process is killed
    instantaneously at 0, and still with $\kappa$ identically 0, they
    claim to show that they can dispense with the assumption
    \eqref{E:extracond1}.  Further results for this case are developed in \cite{MSM01}.
    
    We wish to extend these results to nearly general killing, and
    more general boundary conditions.  For the first part, the
    difficulties are merely technical.  Mandl's proofs rely on W.
    Feller's careful analytical treatment of the diffusion semigroup
    \cite{wF52}.  The relevant parts of Feller's work assume
    $\kappa\equiv 0$, and the corresponding results of It\^o and
    McKean \cite{IM65}, which do allow for general killing, are not
    sufficiently elaborated to be applied in this context.  One part
    of our work will thus be to extend Feller's results rigorously to
    the case of general killing, and then to show that Mandl's methods
    still work in this setting.  This also allows us to dispense with
    Mandl's requirement that the inaccessible boundary be natural; at
    the same time, we are forced to impose the minor additional
    condition (LP).  Our first result is then
\begin{Thm}  \label{T:mainnatreg} 
  If condition (LP) holds, and the initial distribution satisfies
  (ID'), then the process $X_{t}$ converges to the quasistationary
  distribution $\plm$ on compacta.
\end{Thm}

These technical extensions also allow the case of two regular
boundaries with internal killing to be handled.  We define
$\ulm^{(r)}$ as the
$$
\ulm^{(r)}=\inf\ls \gl \, : \, \phi_{\gl} \text{ satisfies
  condition } \eqref{E:fellerboundr}\rs.
$$

\begin{Thm}  \label{T:main2reg}
Suppose that $0$ and $r<\infty$ are both regular boundaries. 
If both boundaries are reflecting and $\kappa\equiv 0$ 
almost everywhere then $\ulm=0$; otherwise, $\ulm$ is positive.
For any Borel subset $A\subset (0,r)$, and any starting
distribution $\mu$,
\begin{equation}
    \ltf e^{\ulm t} \mathbb{P}_\mu \ls X_t\in A\rs     \frac{\int e^{-B(x)}\phi_{\ulm}(x)d\mu(x)\cdot \int_{A}\phi_{\ulm}(y) dy}
{\int_{0}^{r}e^{-B(y)}\phi_{\ulm}^{2}(y)dy}.
    \label{eq:firstlim}
\end{equation}
\end{Thm}

We note that the assumption that $b$ is continuous at the boundary could be dispensed with in the case of two regular endpoints.  In this case, though, we can no longer appeal to general principles (the Krein-Rutman Theorem) to show that there is a positive eigenfunction at the bottom of the spectrum.  We would then need to impose this as a condition:
\begin{equation}  \label{E:extracond}
   \phi_{\ulm^{(r)}} \text{ is nonnegative}.
   \end{equation}

When the interval is infinite, the second part of the problem, proving convergence without conditioning on a compact set, requires a new approach.  We split this part again in two.  One part, which generalizes the results of \cite{CMSM95}, is
\begin{Thm}  \label{T:limsups}
  Suppose (LP) holds, and either
   $$
   \liminf_{z\to\infty} \kappa(z)>\ulm \text{ or }    \limsup_{z\to\infty} \kappa(z)<\ulm.
   $$
Then  either $X_{t}$ converges to the quasistationary distribution $\plm$ from every compactly supported initial distribution, or $X_{t}$ escapes to infinity from every compactly supported initial distribution.  In the case $\liminf_{z\to\infty} \kappa(z)>\ulm$, $X_{t}$ converges to the quasistationary distribution $\plm$ if and only if $\izf \plm(y)dy$ is finite.
\end{Thm}

Note that we have already shown, in Proposition \ref{P:esc}, that the process must escape to infinity if $\izf\plm(y)dy=\infty$.  When the integral is finite, though, there may be ambiguity.

We remark here, as well, that Lemma \ref{L:arelate} tells us that the process has asymptotic killing rate $\ulm$ when it converges to $\plm$.  If $K:=\lzf \kappa(z)$ exists and is finite, the process has asymptotic killing rate $K$ when $X_{t}$ escapes to infinity.  (Unless otherwise indicated, when we say that $\lzf \kappa(z)$ exists, we always mean to include the possibility that $\lzf \kappa(z)=\infty$.)

The general idea here is that the process cannot split, with some mass
going off to infinity and other mass staying near the origin, if the
killing rates differ.  If the killing rate converges to a limit which is larger than $\ulm$, we can be sure which side of the dichotomy we are on, namely, convergence to the quasistationary distribution.

\begin{Thm}  \label{T:limits}
  Suppose the endpoint at $\infty$ is natural, conditions (LP) and
  (GB) hold, and $K:=\lim_{z\to\infty}\kappa(z)$ exists.
\begin{enumerate}  
\item If $K\ne \ulm$, then either
\begin{itemize}
\item  $X_{t}$ has an asymptotic killing rate $\eta=\ulm$ and converges to the quasistationary distribution $\plm$ for every compactly supported initial distribution; or
\item $X_{t}$ has an asymptotic killing rate $\eta=K$ and escapes to infinity for every compactly supported initial distribution.  
\end{itemize}

\item If $K>\ulm$, then $\eta=\ulm$.

\item If $K=\ulm$, then $X_{t}$ has an asymptotic killing rate $\eta=\ulm=K$.

\item Regardless of the relation between $K$ and $\ulm$, we have
\begin{equation}  \label{E:omegalim}
\ltf \omega_{t}(\nu)=\lim_{t\to\infty} \frac{\P_{\nu}\{\tp >t \}}{\P_{1}\{\tp >t \}} = \frac{\psi_{\eta}(x)}{\psi_{\eta}(1)}.
\end{equation}
\end{enumerate}
\end{Thm}

Note that this result dispenses with the constraint (ID') on the initial condition.

In section \ref{sec:example} we present, as an example, a continuous
analogue of the Le Bras process.  The Le Bras process, which was mentioned in section \ref{sec:intro}, is simply a Markov process on the positive integers, in which the transitions are only single steps up, with the transition rate proportional to the state.  The process is killed at a rate which is also proportional to the current state.

For the sake of clarity of exposition, we have placed all lemmas whose
proofs are merely technical, and do not shed light on the central
issues, in section \ref{sec:tech} at the end.

\section{Convergence on compacta}  \label{sec:compact}
\subsection{Proof of Theorem \ref{T:mainnatreg}} \label{sec:mainnatreg}
This result generalizes Theorem 2 of \cite{pM61}. 
The proof there assumed $\kappa\equiv 0$.  We follow the general lines of Mandl's proof, substituting stronger (or clearer) arguments as needed.

\nc{\ho}{\hat{\Omega}} \nc{\tO}{\wt{\Omega}} 

We begin by defining
three differential operators, $\Omega^{*}$, $\ho$, and $\tO$.  All
three act in the same way, as $\mL^{*}$, on different subsets of
$\mC$, the set of continuously twice differentiable functions on
$\R^{+}$ satisfying the boundary condition \eqref{E:fellerbound} at 0.
The subsets are assigned different Banach-space structures.  We use
the following notation: \\[3mm]
\hspace*{-1cm}
\begin{tabular}{c c c c}
Generator & Space & Norm & Domain\\
$\Omega^{*}$ & $L^{1}$ &  $\|f\|_{1}=\izf |f(z)|dz$ & $\mD=\ls f\in
\mC\cap L^{1}\, :\, \Omega f \in L^{1}\rs$\\ 
$\ho$ & $L^{2}(\beta)$ & $\|f\|_{2}=\left(\izf
  |f(z)|^{2}\beta(dz)\right)^{1/2}$ &  $\hat{\mD}=\ls f\in \mC\cap
L^{2}(\beta)\, :\, \ho f \in L^{2}(\beta)\rs$\\ 
$\tO$ & $\ell=L^{1}\cap L^{2}(\beta)$ & $\|f\|_{\ell}\|f\|_{1}+\|f\|_{2}$ & $\wt{\mD}=\ls f\in \ell\, :\, \tO f \in
\ell\rs$ 
\end{tabular}\\[3mm]
We note that all three domains are dense in $\mC$: any function in $\mC$ may be approximated in any of these three norms by a function
that vanishes, along with all its derivatives, outside of a compact
interval.  Mandl's Theorem 1, which we need to prove, states that all
three of these operators generate strongly continuous contraction
semigroups, and that the restrictions of all three to $\ell$ are
identical.

The operator $\Omega^{*}$ is precisely the one analyzed in Lemma
\ref{L:gen}.  There we showed that it generates the adjoint semigroup
of our killed diffusion on $L^1$.

The content of Problem 9.12 of \cite{CL55} is that for any $\ga$ not
in the spectrum of $\ho$ (in particular, for any positive $\ga$), for
$f\in L^{2}(\beta)$,
\begin{equation}  \label{E:psi}
(\ga-\ho)^{-1}f = \Psi^{-1} \left( \frac{\Psi f (\gl)}{\ga+\gl}\right).
\end{equation}
(Note again that the sign on $\gl$ is reversed, because we have
flipped $\varrho$ onto the positive half-line; in particular, we have
defined $\phi_{\gl}$ to have eigenvalue $-\gl$.)  Then, since $\Psi$
is unitary,
$$
\ga\bigl\|(\ga-\ho)^{-1}f\bigr\| = \left\|\frac{\ga}{\ga+\gl} \Psi
  f \right\|\le \| \Psi f\|=\|f\|.
$$
(The inequality depends on the fact that the spectrum of $\ho$
includes no positive values.)  This means that $\ho$ is dissipative,
proving (by the Hille-Yosida Theorem) the existence of a strongly continuous contraction semigroup
on $L^{2}(\beta)$, whose generator is an extension of $\ho$.  To prove
that the generator is actually the closure of $\ho$, we need to show
(by Theorem 2.12 of \cite{EK86}) that the image of $\hat{\mD}$ under
$\ga-\ho$ is dense in $\hat{\mD}$.  Let $S$ be the subset of
$L^{2}(\varrho)$ consisting of functions with compact support.
Consider $(\ga-\ho)^{-1}\Psi^{-1}(S)$.  All functions in this image
are bounded integrals of $\psi(x,\gl)$ (with respect to $\gl$), hence
they are twice differentiable, and the continuity of the diffusion coefficients
implies that they are continuously twice differentiable.  Since
$\Psi^{-1}(S)$ is dense in $L^{2}(\beta)$, this proves that $\mD$ is a
core.

Now choose $f\in \ell$, and any $\gl>0$.  Then (by Proposition 1.2.1 of \cite{EK86}) there exist unique solutions $u_{1}\in L^{1}$ and $u_{2}\in L^{2}(\beta)$ to
\begin{equation}  \label{E:u1}
\gl u_{1} - \Omega^{*} u_{1}= f,
\end{equation}
and
\begin{equation}  \label{E:u2}
\gl u_{2} - \ho u_{2}= f,
\end{equation}
with $\gl\|u_{1}\|_{1}\le \|f\|_{1}$ and $\gl\|u_{2}\|_{2}\le
\|f\|_{2}$.  (Here $\|\cdot\|_{1}$ is the $L^{1}$ Lebesgue norm, and
$\|\cdot\|_{2}$ is the $L^{2}$ norm with respect to $\beta$.)  We show
now that $u_{1}=u_{2}$.

It suffices to assume that $f$ is nonnegative, since we may decompose
$f$ into positive and negative parts.  Since $u_{1}=\int_{0}^{r}
G^{*}_{\ga}(x,y) f(y) dt$, and $G^{*}_{\ga}(x,y)\ge 0$, the function
$u_{1}$ is nonnegative. Equation \eqref{E:psi} implies that $u_{2}$ is
nonnegative as well.

As we have already shown in the proof of Lemma
\ref{L:eigenfunction}, there exists a nonnegative solution $u$ to
the equation $(\mL-\kappa) u - \gl u=0$, satisfying the boundary
condition
$$
(1-p_{0})e^{-B(0)}\phi(0) = p_{0}\phi'(0).  
$$
The function $v(x):=e^{B(x)}g(x)$ is then a nonnegative solution to
$(\mL^{*}-\kappa)v -\gl v=0$, satisfying the boundary condition
\eqref{E:fellerbound}.  Since the solutions to \eqref{E:u1} and
\eqref{E:u2} are unique, it follows that $v$ is neither in $L^{1}$ nor
in $L^{2}(\beta)$.  By uniqueness of solutions to ordinary
differential equations, we must have $u_{1}-u_{2}=cv$, for some $c$.
If $c>0$, then $u_{1}\ge cv \ge 0$, so $u_{1}$ would not be in
$L^{1}$.  On the other hand, if $c<0$, then $u_{2}$ would not be in
$L^{2}(\beta)$.  Consequently, $c=0$, and $u_{1}=u_{2}$ is in $\ell$.
Furthermore,
$$
\|u_{i}\|_{\ell}=\|u_{i}\|_{1}+\|u_{i}\|_{2} \le \|f\|_{1}+\|f\|_{2} =\|f\|_{\ell},
$$
so $\tO$ is dissipative.  It follows from the Hille-Yosida Theorem
\cite[Theorem 1.2.6]{EK86} that $\tO$ generates a strongly continuous
contraction semigroup, for which $\ell$ is a core.  At the same time,
we have shown that the resolvents of all three generators are the same
on $\ell$, so the restrictions of the semigroups to $\ell$ must also
be identical.

We note, now, that \eqref{E:psi} implies, by identity of resolvents,
that the semigroup $\hat{P}_{t}$ generated by $\ho$ maps under $\Psi$
onto the semigroup $Q_{t}$ on $L^{2}(\varrho)$ defined by
$Q_{t}g(\gl)=e^{-\gl t}g(\gl)$.  Consequently, for $f\in \ell$,
\begin{equation}  \label{E:qt}
\Psi P^{*}_{t} f = \Psi \hat{P}_{t}f =e^{-\gl t} \Psi f.
\end{equation}

We may assume, without loss of generality, that (ID) is satisfied, rather than (ID').
Let $f$ be the initial density.  For any bounded Borel set $A$,
by \eqref{E:qt} and the unitarity of $\Psi$, we have
\begin{equation}  \label{E:bmset}
\P \ls X_{t}\in A \rs = \int \indic_{A}(x) P^{*}_{t}f(x) dx= 
\int_{\ulm}^{\infty} e^{-\gl t}\xi_{A}(\gl) g(\gl) d\varrho(\gl),
\end{equation}
where $g:=\Psi f$ and 
$$
\xi_{A}(\gl):=\Psi (e^{B}\indic_{A})(\gl)=\int_{A} e^{B(x)}\phi_{\gl}(x) dx.
$$
Note that
$\xi_{A}(\gl)$ is an integral over a finite interval, hence is continuous.
For any $h\in L^{1}(\varrho)$ define (following Mandl)
$$
\mI(t,h):=\int_{\ulm}^{\infty} e^{-\gl t}h(\gl) d\varrho(\gl).
$$
Following \eqref{E:bmset}, it suffices to show, for bounded
Borel sets $A$ and $A'$,
\begin{equation}  \label{E:sufftoshowmandl}
\ltf \frac{\mI(t,g\xi_{A})}{\mI(t,g\xi_{A'})}\frac{\xi_{A}(\ulm)}{\xi_{A'}(\ulm)}. 
\end{equation}

For $h:(\ulm,\infty)\to\R^{+}$ and $\gl^{*}>\ulm$, define
$\mI(t,h,\gl_{*}):=\mI(t,h\indic_{[\ulm,\gl_{*}]})$.  
We note first that for any $h_1,h_2\in L^1(\varrho)$, with $h_2$
positive, and $\infty\ge \gl_1,\gl'_1,\gl_2>\ulm$,
\begin{align*}
\left|  \frac{\mI(t,h_1,\gl_1)}{\mI(t,h_2,\gl_2)}-
\frac{\mI(t,h_1,\gl'_1)}{\mI(t,h_2,\gl_2)} \right|
\le
   e^{(\ulm-\gl_1\wedge\gl'_1\wedge \gl_2)t}
  \frac{\int_{\ulm}^\infty |h_1(\gl)| d\varrho(\gl)}
   {\int_{\ulm}^\infty h_2(\gl) d\varrho(\gl)},
\end{align*}
so that
$$
\lstf \frac{\mI(t,h_1,\gl_1)}{\mI(t,h_2,\gl_2)} \text{ and }
\litf \frac{\mI(t,h_1,\gl_1)}{\mI(t,h_2,\gl_2)}
$$
are independent of $\gl_1$ and $\gl_2$.  Moreover, we have the bounds
\begin{equation}
  \label{E:lambdacomp}
\limsup_{\gl\downarrow\ulm} \frac{h_1(\gl)}{h_2(\gl)} 
\ge \lstf \frac{\mI(t,h_1,\gl_1)}{\mI(t,h_2,\gl_2)} \ge
\litf \frac{\mI(t,h_1,\gl_1)}{\mI(t,h_2,\gl_2)}\ge
\liminf_{\gl\downarrow\ulm} \frac{h_1(\gl)}{h_2(\gl)} 
\end{equation}
In particular, assumption (ID) implies
\begin{equation}
  \label{E:idlim}
  \sup_{\gl_*>\ulm} \lstf \frac{\mI(t,g^-)}{\mI(t,1,\gl_*)}<\infty,
\end{equation}
where $g^{-}=\max\{-g,0\}$ and $g^{+}=\max\{g,0\}$.

Suppose that $f$ has compact support.  Then $g$ is continuous, and
$$
g(\ulm) = \int_{0}^{\infty}f(x) \phi_{\ulm}(x) dx>0.
$$
Thus, for any finite $\gl_{*}>\ulm$,
$$
\ltf \frac{\mI(t,g\xi_{A})}{\mI(t,1,\gl_{*})} = \xi_{A}(\ulm)g(\ulm)>0,
$$
from which \eqref{E:sufftoshowmandl} is immediate.

Now, consider a general density $f$.  Choose $z_{0}$ such that $\int_{0}^{z_{0}}f(y) dy>0$, let $f_{0}=f\indic_{[0,z_{0}]}$, and let $g_{0}=\Psi f_{0}$.  Then
\begin{equation}  \label{E:compcomp}
\liminf_{t\to\infty} \frac{\mI(t,g\xi_{A})}{\mI(t,1,\lambda_{*})}\ge \ltf \frac{\mI(t,g_{0}\xi_{A})}{\mI(t,1,\gl_{*})}  >0.
\end{equation}

We claim that
\begin{equation}
  \label{E:idgd}
  \gd:= 1- \sup_{\gl_*>\ulm} \lstf \frac{\mI(t,g^-)}{\mI(t,g^+,\gl_*)}>0.
\end{equation}
Suppose this is true.  Then
\begin{equation}
  \label{eq:idgd2}
  \begin{split}
  \lstf \frac{\mI(t,|g|)}{\mI(t,g,\gl_*)}&\le 
  \lstf \left(\frac{\mI(t,g^{+},\gl_{*})-\mI(t,g^{-},\gl_{*})}{\mI(t,g^{+})+\mI(t,g^{-})}\right)^{-1}\\
  &\le \left(1- \lstf \frac{\mI(t,g^-,\gl_{*})}{\mI(t,g^+)} \right)^{-1}\\
  &\le\gd^{-1}.
  \end{split}
\end{equation}
We may then conclude that for any continuous function $h:[\ulm,\infty)\to\R$,
\begin{align*}
  h(\ulm)-\gd^{-1} \sup_{\gl\in [\ulm,\gl_*]}\lv h(\gl) - h(\ulm)\rv &\le
\litf \frac{\mI(t,gh)}{\mI(t,g,\gl_*)}\\
&\le \lstf \frac{\mI(t,gh)}{\mI(t,g,\gl_*)} \\
&\le h(\ulm)+\gd^{-1} \sup_{\gl\in [\ulm,\gl_*]}\lv h(\gl) - h(\ulm)\rv .
\end{align*}
Thus, for any $\gl_*>\ulm$, and continuous $h_{1}$ and $h_{2}$ with  $h_{2}(\ulm)>\gd^{-1}\sup_{\gl\in[\ulm,\gl_*]} |h_{2}(\gl)-h_{2}(\ulm)|$,
\begin{multline*}
\frac{\gd h_{1}(\ulm) - \sup_{\gl\in[\ulm,\gl_*]}
  |h_{1}(\gl)-h_{1}(\ulm)|}
{\gd h_{2}(\ulm) +\sup_{\gl\in[\ulm,\gl_*]}
  |h_{2}(\gl)-h_{2}(\ulm)|}\\
\le \litf \frac{\mI(t,gh_{1})}{\mI(t,gh_{2})}
\le \lstf \frac{\mI(t,gh_{1})}{\mI(t,gh_{2})} \\
\le \frac{\gd h_{1}(\ulm) + \sup_{\gl\in[\ulm,\gl_*]}
  |h_{1}(\gl)-h_{1}(\ulm)|}
{\gd h_{2}(\ulm) -\sup_{\gl\in[\ulm,\gl_*]}
  |h_{2}(\gl)-h_{2}(\ulm)|}.
\end{multline*}
Sending $\gl_*$ to $\ulm$, we see that
\begin{equation}  \label{E:uselater}
\ltf \frac{\mI(t,gh_{1})}{\mI(t,gh_{2})} = \frac{h_{1}(\ulm)}{h_{2}(\ulm)}.
\end{equation}
Applying this with $h_{1}=\xi_{A}$ and $h_{2}=\xi_{A'}$ completes the proof of \eqref{E:bmset}.  

It remains to prove the claim \eqref{E:idgd}.  By continuity of $\xi_A$, and
\eqref{E:idlim}, we may find $\gl^*>\ulm$ and a positive constant
$\gep$ such that for all $\gl_*\in (\ulm,\gl^{*})$,
$$
\litf \frac{\mI(t,g\xi_{A},\lambda_{+})}{\mI(t,1,\gl_{*})}-
2\sup_{\gl\in (\ulm,\gl^{*})}\lv \xi_{A}(\gl) - \xi_{A}(\ulm) \rv
\cdot \sup_{\gl\in [\ulm,\gl_{*}]} g^{-}(\gl) >
\gep \lstf \frac{\mI(t,g^{-})}{\mI(t,1)}.
$$
Define $\delta(\gep):=\gep\xi_{A}(\ulm)/(1+\gep\xi_{A}(\ulm))$.  Suppose, for some
$\gl_{*}<\gl^{*}$, there were a sequence $t_{n}\to \infty$ such that
$$
\lnf  \frac{\mI(t_{n},g^{-})}{\mI(t_{n},g^{+},\gl_{*})}>1-\gd(\gep).
$$
Then
\begin{align}
\lsnf \frac{\mI(t_{n},|g|)}{\mI(t_{n},1,\gl_{*})}&<
   \frac{2}{1-\gd(\gep)}\lsnf \frac{\mI(t_{n},g^{-})}{\mI(t_{n},1,\gl_{*})}\text{ and} \label{E:2gdbd1}\\
\lsnf \frac{\mI(t_{n},g)}{\mI(t_{n},1,\gl_{*})}&<
   \frac{\gd(\gep)}{1-\gd(\gep)}\lsnf \frac{\mI(t_{n},g^{-})}{\mI(t_{n},1,\gl_{*})} .\label{E:2gdbd2}
\end{align}
We may conclude that
\begin{align*}
\linf \frac{\mI(t_{n},g\xi_{A})}{\mI(t_{n},1,\gl_{*})}
&\le \xi_{A}(\ulm) \linf \frac{\mI(t_{n},g)}{\mI(t_{n},1,\gl_{*})}\\
  &\hspace*{1cm} + \sup_{\gl\in (\ulm,\gl^{*})}\lv  \xi_{A}(\gl) - \xi_{A}(\ulm) \rv\cdot
   \frac{2}{1-\gd(\gep)}\lsnf \frac{\mI(t_{n},g^{-})}{\mI(t_{n},1,\gl_{*})},
\end{align*}
so that
$$
\linf \frac{\mI(t_{n},g)}{\mI(t_{n},1,\gl_{*})} \ge 
   \frac{\gep}{\xi_{A}(\ulm)} \lsnf \frac{\mI(t_{n},g^{-})}{\mI(t_{n},1,\gl_{*})},
$$
contradicting \eqref{E:2gdbd2}.  We have thus proven the claim \eqref{E:idgd} for $\gl_{*}$ sufficiently close to $\ulm$.  By monotonicity, this completes the proof.

\subsection{Proof of Theorem \ref{T:main2reg}}\label{sec:main2reg}
\nc{\irr}{\int_{0}^{r}}
Let $G_{\ga}(x,y)$ be the $\ga$-Green's function for the operator $\mL$, with respect to Lebesgue measure.  The Green's function for $\mL^{*}$ is then given by $G^{*}_{\ga}(x,y)=G_{\ga}(y,x)$.  Equation 6 of section 4.11 of \cite{IM65} gives us an explicit representation of the Green's function $G_{\ga}(x,y)$, from which we can derive (by Theorem VI.23 of \cite{RS721}) that the operator on $L^2(\beta)$
generated by the kernel $G_{\ga}^{*}$ is compact.  It is also self-adjoint with 
respect to the measure $\beta$.  By the Hilbert-Schmidt Theorem 
(Theorem VI.16 of \cite{RS721}), the 
spectrum is discrete, and has only 0 as accumulation point.  
Furthermore, there is a complete orthonormal basis of eigenfunctions.  
These are also eigenfunctions for the infinitesimal generator; the 
eigenvalue $-\gl$ for the generator corresponds to $\rec{\mu+\gl}$ for 
$G_{\mu}$.  Let $-\gl_{0}\ge -\gl_{1}\ge\cdots$ be the eigenvalues, 
and $\phi_{k}$ the eigenfunction of the generator \eqref{E:ivp} 
corresponding to $\gl_{k}$. Since the 
eigenfunctions are orthogonal with respect to the positive 
function $e^{-B(x)}$, no more than one may be nonnegative.  Either by 
assumption, or (if the drift is finite at the boundaries)
by the Krein-Rutman Theorem in the form derived as 
Theorem 2.1 of \cite{rN84}, there is a nonnegative eigenfunction, and 
it corresponds to the top eigenvalue $-\gl_0=-\ulm$.   (With a slight 
further restriction of the conditions, we could instead
 apply the more elementary result, Theorem 8.2.1 of \cite{CL55}.)

\nc{\htk}{h_{t}^{K}}
If $h(x)$ is a smooth function on $[0,r]$ with
$$
\|h\|^2:=\irr h(x)^{2} \beta(dx)<\infty,
$$
we can write
$$
h=\sum_{k=0}^{\infty} \left( \irr h(y) \phi_{k}(y)\beta(dy) 
\right) \phi_{k} =:\sum_{k=0}^{\infty} a_{k} \phi_{k}
$$
(The sum need not converge pointwise, but in $L^2(\beta)$.)
Let 
$$
\htk(x)=\sum_{k=0}^{K} e^{-\gl_{k}t}a_{k}\phi_{k}(x).
$$

If the initial distribution has a smooth density $h$, then by the Cauchy-Schwartz inequality (following an equivalent computation on page 561 of \cite{pM61}),
\begin{align*}
    \ltf \Bigl|e^{\ulm t}\Pr X_{t}\in A\rs - 
    a_{0}&\int_{A}\phi_{\ulm}(x) dx\Bigr| \\
    &    \ltf \lim_{K\to\infty} \int_{A} e^{-B(x)/2}\lv e^{\ulm t}\htk - 
    a_{0}\phi_{0} \rv e^{B(x)/2} dx\\
    &\le \ltf \lim_{K\to\infty} \left\|e^{\ulm t}\htk - a_{0}\phi_{0}\right\|
    \cdot \irr e^{B(x)}dx\\
    &=0.
\end{align*}

If the initial distribution does not have a smooth density, we start instead with the distribution at a small finite time $\gep$.

\section{Convergence on $\R^{+}$} \label{sec:infconv}
\subsection{The limit functions $\omega$ and $F$}  \label{sec:lfo}
\nc{\bom}{\bar{\omega}}
\begin{Lem}\label{L:omegat}
Given any increasing sequence $t^{*}_{n}\to\infty$, we may find a subsequence $t_{n}\to\infty$ with the following properties:
\begin{enumerate}
\item There is a limit $\omega^{(t_{n})}(x)=\lnf \omega_{t_{n}}(x)$ in $C(\R^{+})$, with convergence in the topology of uniform convergence on compacta.  The function $\omega$ is nonzero everywhere, except when the process is killed instantaneously at 0, in which case $\omega(0)=0$.

\item There is a limit $F^{(t_{n})}=\lnf F_{t_{n}}$ in $C(\mc,\mo)$, supplied with the topology of uniform convergence on compact sets.
\end{enumerate}
When there is no need to emphasize the particular subsequence, we will drop the superscript $(t_{n})$ from the notation.
\end{Lem}

\begin{proof}
We consider first the functions $F_{t}$ and $\omega_{t}$ restricted to $\delta$-measures.  The topology on $\mc$ restricts then to the usual topology on $\R^{+}$.

For integers $n\ge 1$, $t\in [0,\infty)$, and $x\in [n,n+1]$, define
\begin{equation}  \label{E:omegan}
\otn(x):=\frac{\P_{x} \{\tp>t\}}{\P_{n} \{\tp>t\}}.
\end{equation}
We show that for each $n$ the set of functions $\otn$ is uniformly equicontinuous.  The bound \eqref{E:unifomegabd} immediately tells us that
$$
\otn(x)\le \P_{n}\ls \tau_{n+1}<\tp \rs^{-1}<\infty
$$
for all $x$, so the functions are uniformly bounded.  If $x,y\in [n,n+1]$, then
$$
\lv \otn(x) -\otn(y)\rv \le \P_{n}\ls \tau_{n+1}<\tp \rs^{-1} \Bigl[ \P_{x} \ls \tau_{y}<\tp \rs^{-1}\vee 
\P_{y} \ls \tau_{x}<\tp \rs^{-1} \, -1\Bigr].
$$
Since the diffusion rate is constant, and the drift and killing are bounded on $[n,n+1]$, this bound goes to zero uniformly as $|x-y|\to 0$.

By Ascoli's Theorem \cite[Theorem 0.4.13]{rE95}, the set of functions $\{\otn\}_{t\in\R^{+}}$ is relatively compact.  Thus, any sequence $t'_{i}\to\infty$ has a subsequence $(t_{i})$, such that $\phantom{}_{n}\omega_{t_{i}}$ converges uniformly to a function $\omega^{(n)}:[n,n+1]\to \R^{+}$.  Furthermore, we may take the sequence for $\otn$ to be a refinement of that for 
$\phantom{}_{n-1}\omega_{t}$.  By taking a diagonal subsequence, we then get a single sequence $t_{i}$ such that $\phantom{}_{n}\omega_{t_{i}}$ converges uniformly to $\phantom{}_{n}\omega$ for every $n$.  By an identical argument, we can also assume that $\phantom{}_{0}\omega_{t}:=\P_{x} \{\tp>t\}/\P_{1} \{\tp>t\}$ for $x\in [0,1]$ also converges along the sequence $t_{i}$.  (We treat the case $x<1$ separately only because killing at 0 could make it inconvenient to place $\P_{0}\{\tp>t\}$ in the denominator.)  If we now define $\omega_{t}$ by
$\omega_{t}(x) = \P_{x} \{\tp>t\}/\P_{1} \{\tp>t\}$ for all $x\in [0,\infty)$, we see that $\omega_{t}$ is a product of finitely many functions $\otn$, so that $\omega=\lim_{i\to\infty} \omega_{t_{i}}$ exists as a positive continuous function, though the convergence need only be uniform on compacta.

We note, as well, that for all $y\in\R^{+}$, by \eqref{E:unifomegabd},
$$
\omega(y)\ge \P_{y}\ls \tau_{1}<\tp\rs,
$$
which is positive, except when $y=0$ and there is complete killing at 0.

\nc{\nf}{\phantom{}_{n}F}
We continue with $F_{t}$.  We note that it will suffice to prove (ii) for the measures $\nf_{t}$, defined by
$$
\nf_{t}(x,A)=F_{t}(x,A\cap [\rec{n},n]);
$$
the measure $F(x,\cdot)$ is then defined as the limit along a common refinement of the subsequences defining $\nf(x,\cdot)$.  We need to show that the set of functions $F_{t}$ is relatively compact.  Since relative compactness is a topological property, it will suffice to prove equicontinuity in any convenient metric that metrizes the compact weak-* topology.  We choose the Wasserstein metric \cite{sR92}, given by
$$
d(\mu,\mu')= \sup \ls \lv \mu(g)-\mu'(g)\rv \rs,
$$
where the supremum is taken over functions $g:\R^{+} \to [0,1]$, supported on $[\rec{n},n]$, and with Lipschitz constant $\le 1$.

Since the space is bounded, we need only to show that for any $x\in [0,\infty)$ (or $x\in (0,\infty)$, in the case of complete killing at 0),
\begin{equation}  \label{E:equix}
\lim_{x'\to x} \sup_t d\lp \nf_{t}(x,\cdot) , \nf_{t}(x',\cdot)\rp =0.
\end{equation}
Let $R_{x',x}$ be the (subprobability) distribution of $\tau_{x}$, when the process is started from $x'$.  Let $g$ be any appropriate test function.  

By the strong Markov property, for any $t_{0}$,
\begin{align*}
 \Bigl| \int \nf_{t}(x&,dy)g(y) - \int \nf_{t}(x',dy)g(y) \Bigr|\\
&\le \frac{\int_{0}^{t_{0}\wedge t} R_{x',x}(ds) \lv\E_{x} g(X_{t-s}) - E_{x} g(X_{t})\rv}{\P_{x}\{ \tp>t\}}
+ \left| \frac{\P_{x'}\{\tp>t\}}{\P_{x}\{\tp>t\}} - 1 \right|\\
&\hspace*{.5cm}+R_{x',x}\lp [t_{0}\wedge t,\infty]\rp 
    \biggl|\E_{x'}\lb g(X_{t})\cond \tp>t \, \& \, \tau_{x}>t_{0}\wedge t \rb - 
    \E_{x}\lb g(X_{t}) \cond \tau>t\rb \biggr|
 \end{align*}
 \begin{multline*}
\le \omega_{**}(x,x')+ \int_{0}^{t_{0}} R_{x',x}(ds) \P_{x} \ls X_{t-s}> n+1 \text{ and } X_{t}\le n \cond \tp>t \rs\\
+ \int_{0}^{t_{0}} R_{x',x}(ds) \E_{x} \Bigl[  (1\wedge|X_{t}-X_{t-s}|) \indic\ls X_{t-s}\le n+1\rs \cond \tp>t \Bigr]\\
+R_{x',x}\lp [t_{0}\wedge t,\infty]\rp \biggl(
    \E_{x'}\lb 1\wedge\lv X_{t} - x' |\cond  \tp>t \, \& \, \tau_{x}>t_{0}\wedge t  \rb \\
+ \E_{x} \lb 1\wedge |X_{t}-x| \rb +|x-x'|\biggr) .
\end{multline*}    
Note that this bound is independent of $g$.  It remains, then, only to show that the supremum over $t$ converges to 0 as $x'\to x$.  By regularity, $R_{x',x}$ converges to $\delta_{0}$.    Since $\omega_{**}$ is continuous, and $\omega_{**}(x,x)=0$, the first term goes to 0 as $x'\to x$.   The second term is bounded by
\begin{equation}  \label{E:firstbound}
\int_{0}^{t} R_{\nu}(ds) \sup_{y\ge n+1}\P_{y} \ls X_{s}\le n \cond \tp>s \rs.
\end{equation}
The third term is bounded by
\begin{equation}  \label{E:secondbound}
\int_{0}^{t} R_{\nu}(ds) \sup_{y\le n+1}\E_{y} \lb  (1\wedge|X_{s}-y|)  \sup_{y\le n+1} \P_{y}\ls \tp>s \rs^{-1}.
\end{equation}
The Bounded Convergence Theorem implies then --- since $X$ is a Feller process --- that \eqref{E:firstbound} and \eqref{E:secondbound} both converge to 0 as $R_{\nu}\to\delta_{0}$.

We may conclude, for any positive $t_{0}$, that
\begin{align*}
\lim_{x'\to x}\sup_{t}  &\Bigl| \int \nf_{t}(x,dy)g(y) - \int \nf_{t}(x',dy)g(y) \Bigr|\\
&\le \lim_{x'\to x}\sup_{t} R_{x',x}\lp [t_{0}\wedge t,\infty]\rp\\
&\hspace*{5mm}\times \Bigl(
    \E_{x'}\lb 1\wedge\lv X_{t} - x' |\cond \tp\wedge \tau_{\nu}>t \rb 
    + \E_{x} \lb 1\wedge |X_{t}-x| \rb +|x-x'|\Bigr)\\
&\le \lim_{x'\to x} 2 R_{x',x}\lp [t_{0},\infty]\rp +  
\lim_{x'\to x}\sup_{t\le t_{0}} R_{x',x}\lp [t,\infty]\rp\E_{x} \lb 1\wedge |X_{t}-x| \rb\\
&\hspace*{1cm}   +  \lim_{x'\to x}\sup_{t\le t_{0}}R_{x',x}\lp [t,\infty]\rp\E_{x'}\lb 1\wedge\lv X_{t} - x' |\cond  \tp>t \, \& \, \tau_{x}>t_{0}\wedge t \rb .
\end{align*}
By continuity of the diffusion parameters, we may find an increasing function $h(t)\in [0,1]$ such that $\lim_{t\to 0} h(t)=0$, and for all $x'\le x+1$,
\begin{align*}
\E_{x'}\lb 1\wedge\lv X_{t} - x' |\cond  \tp>t \, \& \, \tau_{x}>t_{0}\wedge t  \rb &\le h(t) \text{ and}\\
\E_{x}\lb 1\wedge\lv X_{t} - x |\cond \tp>t\rb &\le h(t).
\end{align*}
The bound may then be turned into
$$
\lim_{x'\to x}\sup_{t\le t_{0}} 2R_{x',x}\lp [t,\infty]\rp h(t) \le 2h(t_{1})+2\lim_{x'\to x} R_{x',x}\lp [t_{1},\infty]\rp \le 2h(t_{1})
$$
for any $t_{1}>0$, which goes to 0 as $t_{1}\to 0$.

It remains to extend the result to all $\nu\in\mc$.  We define $\omega$ and $F$ by the formulae \eqref{E:omeganu} and \eqref{E:Fnu}.  Since $\omega(x)>0$ for $x>0$, $\omega(\nu)$ is nonzero, so that $F(\nu,\cdot)$ is well defined.  Since $\omega$ is bounded on compact sets, for any sequence $t_{n}$ such that $\lnf \omega_{t_{n}}(x)=\omega(x)$, we also have $\lnf \omega_{t_{n}}(\nu)=\omega(\nu)$, by the Bounded Convergence Theorem.  If $\nu_{n}$ is a sequence of measures converging to $\nu$ in $\mc$, they converge weakly and have their supports contained in a  common compact set.  Since $\omega$ is continuous, it follows again from the definition of weak convergence that $\lnf\int \omega(x)d\nu_{n}(x) =\int \omega(x)d\nu(x)$.

Now, consider any sequence $t_{n}\to\infty$, such that $\lnf F_{t_{n}}(x,\cdot)=F(x,\cdot)$.  The same bounded convergence argument shows that $\lnf F_{t_{n}}(\nu,g)=F(\nu,g)$ for any $\nu\in\mc$ and any continuous $g:\R^{+}\to [0,1]$.  Suppose $\nu_{n}$ is a sequence of measures converging to $\nu$ in $\mc$.  Then
$$
\lnf \omega(\nu_{n})F(\nu_{n},g) = \lnf \izf \omega(x) F(x,g) d\nu_{n}(x)=\omega(\nu)F(\nu,g),
$$
since $x\mapsto\omega(x) F(x,g)$ is continuous and bounded on the compact set containing the supports of the measures $\nu_{n}$.  Since $\omega(\nu)$ is nonzero, it follows that $F$ is continuous.
\end{proof}

\nc{\oK}{\bar{K}}

\subsection{Invariance properties}  \label{sec:conres}
\begin{Lem}
\label{L:whichdensity}
Suppose $X_{t}$ converges on compacta to the quasistationary distribution $\phi$, when started at the compactly supported distribution $\nu$.  Let $t_{n}$ be any sequence such that $F_{t_{n}}$ has a continuous limit, as described in Lemma \ref{L:omegat}.  Then $F^{(t_{n})}(\nu,\cdot)$ has density $F^{(t_{n})}(\nu,\R^{+}) \phi$ with respect to Lebesgue measure.  If convergence on compacta holds for a set of initial distributions which is dense in $\mc$ (in particular, the distributions satisfying (IDC)), then $F^{(t_{n})}(\nu,\cdot)$ has density $F^{(t_{n})}(\nu,\R^{+}) \phi$ for all $\nu\in\mc$.
\end{Lem}

\begin{proof}
Convergence on compacta implies that for any $z,z'>0$,
$$
F^{(t_{n})}(\nu,[0,z])=F^{(t_{n})}(\nu,[0,z']) \frac{\int_{0}^{z} \phi(y)dy}{\int_{0}^{z'} \phi(y)dy}.
$$
We may normalize $\phi$ so that $\izf \plm(y)dy$ is either 1 or $\infty$.  Letting $z'\to\infty$, we get
\begin{equation}  \label{E:bothsidescont}
F^{(t_{n})}(\nu,[0,z])=F^{(t_{n})}(\nu,\R^{+}) \int_{0}^{z} \plm(y)dy
\end{equation}
for any $\nu$ with a compactly supported density, and any appropriate sequence $t_{n}$.  By continuity, if we take a weakly convergent sequence of (IDC) distributions $\nu$, we see that \eqref{E:bothsidescont} holds for all compactly supported $\nu$.
\end{proof}

\begin{Lem}  \label{L:arelate}
Suppose $X_{t}$ converges on compacta to the quasistationary distribution $\phi$, when started at the compactly supported distribution $\nu$. 
If $F^{(t_{n})}(\nu,\R^{+})=1$, then the limit $\lnf a_{t_{n}}(\nu,r)$ exists, and is equal to $\izf \phi(y)\P_{y}\ls \tp>r\rs dy$.  More generally, if $K:=\lzf \kappa(z)$ exists, then the limit $\lnf a_{t_{n}}(\nu,r)$ exists, regardless of the value of $F^{(t_{n})}(\nu,\R^{+})$ (recall the notation from \eqref{E:definea}), and is equal to 
\begin{equation}  \label{E:arelate}
a^{(t_{n})}(\nu,r)=F^{(t_{n})}(\nu,\R^{+})\int \phi(y)\P_{y}\ls \tp>r\rs dy +\lp 1- F^{(t_{n})}(\nu,\R^{+}) \rp e^{-Kr}.
\end{equation}
\end{Lem}

\begin{proof}
By Lemma \ref{L:whichdensity}, $F^{(t_{n})}(\nu,\cdot)$ has density $F^{(t_{n})}(\nu,\R^{+}) \plm$ with respect to Lebesgue measure.  If $g:\R^{+}\to\R$ is any bounded continuous function such that $g(\infty)=\lim_{y\to\infty}g(y)$ exists, weak convergence implies
\begin{equation}  \label{E:wc}
\lnf \izf F_{t_{n}}(\nu,dy)g(y) = \izf  F^{(t_{n})}(\nu,dy)g(y) + 
   \lp 1-F^{(t_{n})}(\nu,\R^{+}) \rp g(\infty).
\end{equation}
In general,
$$
a^{(t_{n})}(\nu,r)=\int F^{(t_{n})}(\nu,dy) \P_{y}\ls\tp>r\rs.
$$
If $\kappa$ has a limit at $\infty$, the function $y\mapsto \P_{y} \{ \tp>r \}$ is continuous at $\infty$ for each fixed $r$, from which \eqref{E:arelate} follows.  If $F^{(t_{n})}(\nu,\R^{+})=1$, then the term at infinity vanishes in any case, so the convergence of $\kappa$ is irrelevant.
\end{proof}

\begin{Prop}  \label{P:harmonic}
Suppose the conditions (LP) and (GB), are satisfied.  Let $(t_{n})$ be any sequence such that $F_{t_{n}}$ and $\omega_{t_{n}}$ converge to continuous limits, as described in Lemma \ref{L:omegat}.  Then 
\begin{enumerate}
\item  The product $\omega^{(t_{n})}(x)F^{(t_{n})}(x,\R^{+})$ is $\ulm$-invariant.  It satisfies
\begin{equation}  \label{E:qom1}
\omega^{(t_{n})}(x) F^{(t_{n})}(x,\R^{+}) = F^{(t_{n})}(1,\R^{+})\frac{\psi_{\ulm}(x)}{\psi_{\ulm}(1)}.
\end{equation}

\item If $K:=\lzf \kappa(z)$ exists, then $\omega^{(t_{n})}(x)(1-F^{(t_{n})}(x,\R^{+}))$ is $K$-invariant.  If $K$ is finite, then 
\begin{equation}  \label{E:qom2}
\omega^{(t_{n})}(x) \lp 1-F^{(t_{n})}(x,\R^{+})\rp = \lp 1- F^{(t_{n})}(1,\R^{+}) \rp \frac{\psi_{K}(x)}{\psi_{K}(1)}.
\end{equation}
If $K=\infty$ then $F^{(t_{n})}(x,\R^{+})=1$ for all $x$.

\item If $K:=\lzf \kappa(z)$ exists and is finite, then $\omega^{(t_{n})}$ is a convex combination of $\psi_{\ulm}$ and $\psi_{K}$.  Specifically,
\begin{equation}  \label{E:omegaform}
\omega^{(t_{n})}(x)= F^{(t_{n})}(1,\R^{+}) \frac{\psi_{\ulm}(x)}{\psi_{\ulm}(1)}+(1-F^{(t_{n})}(1,\R^{+})) \frac{\psi_{K}(x)}{\psi_{K}(1)}.
\end{equation}

\item If $K:=\lzf \kappa(z)$ exists, then for any compactly supported $\nu$ and $s>0$,
\begin{equation}  \label{E:qom3}
\E_{\nu}\omega^{(t_{n})}(X_{s}) = a^{(t_{n})}(\nu,s) \omega^{(t_{n})}(\nu).
\end{equation}

\item Either $X_{t}$ escapes to infinity from every compactly supported initial distribution, or $\plm$ is the density of a $\ulm$-invariant distribution for $X_{t}$; that is, 
\begin{equation}  \label{E:lambdaharm}
P^{*}_{t}\plm=e^{-\ulm t}\plm.
\end{equation}
\end{enumerate}
\end{Prop}

\begin{proof}
Condition (GB) allows us to apply the Bounded Convergence Theorem, to show that for any positive $s$,
\begin{equation}  \label{E:qinvar0}
\begin{split}
\E_{\nu}\omega^{(t_{n})}(X_{s})&=\E_{\nu}\lnf \frac{\P_{X_{s}} \{\tp>t_{n}\}}{\P_{1}\{\tp>t_{n}\}}\\
&=\omega^{(t_{n})}(\nu)\lnf \P_{\nu} \{\tp>t_{n}+s\cond\tp>t_{n}\},
\end{split}
\end{equation}
which proves \eqref{E:qom3}.

Define
$$
\gm(x):=\omega^{(t_{n})}(x)F^{(t_{n})}(x,\R^{+}).
$$
Then
\begin{equation}  \label{E:qinvar}
\begin{split}
\E_{x} \gm(X_{s}) &= \E_{x} \lzf \lnf \frac{\P_{X_{s}} \{X_{t_{n}}\le z \}}{\P_{1}\{\tp>t_{n}\}}\\
&=\omega^{(t_{n})}(x)\lzf \lnf \frac{\P_x \{X_{t_{n}+s}\le z \}}{\P_{x}\{\tp>t_{n}\}}\\
&=\omega^{(t_{n})}(x)\lzf \lnf \izf F_{t_{n}}(x,dy) \P_{y} \ls X_{s}\le z \rs\\
&= \omega^{(t_{n})}(x)\lzf \izf F(x,dy) \P_{y}\ls X_{s} \le z \rs\\
&=\omega^{(t_{n})}(x)F(x,\R^{+}) \izf \plm(y)\P_{y}\ls \tp>s \rs dy\\
&=\gm(x)\izf \plm(y)\P_{y}\ls \tp>s \rs dy.
\end{split}
\end{equation}
Note that the step from the third to the fourth lines uses the fact that the function $y\mapsto\P_{y} \ls X_{s}\le z \rs$ vanishes at $\infty$ for each $z$, while the following step uses Lemma \ref{L:whichdensity}.

Suppose there is some initial distribution $\nu$ satisfying (IDC), from which $X_{t}$ does not escape to $\infty$.  Then Proposition \ref{P:esc} tells us that $\izf \plm (z)dz$ is finite.  There is some sequence $t_{n}$ such that $F^{(t_{n})}(\nu,\R^{+})>0$.  By the relation \eqref{E:Fnu}, it follows that there must be some $x$ such that $F_{(t_{n})}(x,\R^{+})>0$.  Given any positive $s$ and $t$,
\begin{align*}
\E_{x} \gm(X_{s+t})&= \E_{x}\E_{X_{t}} \gm(X_{s}) \\
&=\E_{x} \gm(X'_{t})\izf \plm(y)\P_{y}\ls \tp>s \rs dy\\
&=\gm(x)\izf \plm(y)\P_{y}\ls \tp>t \rs\izf \plm(y)\P_{y}\ls \tp>s \rs dy.
\end{align*}
Since this is also
$$
\gm(x)\izf \plm(y)\P_{y}\ls \tp>s+t \rs dy,
$$
it follows that the function
$$
s\mapsto \izf \plm(y)\P_{y}\ls \tp>s \rs dy
$$
is multiplicative, and so may be written in the form $e^{-\eta s}$ for some $\eta$.
Observe now that $\gm(x)>0$.  Since $\gm$ is an $\eta$-invariant function, it must be in the domain of $\mL$, and satisfy $\mL\psi = -\eta \psi$.  By Lemma \ref{L:bottom}, and the fact that $\gm$ is nonnegative and not identically 0, it must be that $\eta\le \ulm$.  Lemma \ref{L:subinvar} shows then that $\plm$ is $\ulm$-invariant, proving \eqref{E:lambdaharm}.  This also proves \eqref{E:qom1}, in the case when $X_{t}$ does not escape to infinity from every initial distribution.

Suppose now that $F(\nu,\R^{+})=0$ for all (IDC) distributions $\nu$.  By continuity, $F(x,\R^{+})=0$ for all $x$.  Then \eqref{E:qom1} holds trivially.

\nc{\tg}{\tilde{\gm}}
It remains only to show \eqref{E:qom2} --- equation \eqref{E:omegaform} merely combines this with \eqref{E:qom1}.  We repeat the computation of \eqref{E:qinvar}, now with
$$
\tg(x):=\omega^{(t_{n})}(x)\lp 1-F^{(t_{n})}(x,\R^{+})\rp.
$$
We have
\begin{equation}  \label{E:qinvar2}
\begin{split}
\E_{x} \tg(X_{s}) &= \E_{x} \lzf \lnf \frac{\P_{X_{s}} \{X_{t_{n}}> z \}}{\P_{1}\{\tp>t_{n}\}}\\
&= \omega^{(t_{n})}(x)\biggl[\lzf \int F^{(t_{n})}(x,dy) \P_{y}\ls X_{s} > z \rs \\
&\hspace*{3cm}+ \lp 1- F^{(t_{n})}(x,\R^{+})\rp \lim_{y\to\infty} \P_{y} \ls X_{s}>z\rs\biggr]\\
&=\tg(x)e^{-Ks}.
\end{split}
\end{equation}
In the last step, we have used the fact that $\infty$ is a natural boundary, which implies that
$$
\lim_{y\to\infty} \P_{y} \ls X_{s}>z\rs=\lim_{y\to\infty} \P_{y} \ls \tp>s\rs.
$$
\end{proof}

\begin{Cor}
\label{C:solidarity}
Suppose (LP) and (GB) hold. Let $t_{n}\to\infty$ be a sequence such that $F_{t_{n}}$ and $\omega_{t_{n}}$ have continuous limits, as described in Lemma \ref{L:omegat}. Then either $F^{(t_{n})}(\nu,\R^{+})=0$ for all compactly supported $\nu$ or $F^{(t_{n})}(\nu,\R^{+})>0$ for all compactly supported $\nu$.  If $K=\lzf \kappa(z)$ exists, then either $F^{(t_{n})}(\nu,\R^{+})=1$ for all compactly supported $\nu$ or
$F^{(t_{n})}(\nu,\R^{+})<1$ for all compactly supported $\nu$.
\end{Cor}

\begin{proof}
Suppose there is some $\nu$ such that $F^{(t_{n})}(\nu,\R^{+})=0$.  By \eqref{E:Fnu}, there must be some $x$ such that $F^{(t_{n})}(x,\R^{+})=0$.  By Proposition \ref{P:harmonic}, we know that
$$
\E_{x} \omega^{(t_{n})}(X_{s})F^{(t_{n})}(X_{s},\R^{+})=0.
$$
By regularity, and the fact that $\omega$ is positive, it follows that $F^{(t_{n})}(x,\R^{+})=0$ for Lebesgue-almost-every $x$.  Suppose there were exceptional $x'$ such that $F^{(t_{n})}(x',\R^{+})>0$.  Then
$$
\E_{x'} \omega^{(t_{n})}(X_{s})F^{(t_{n})}(X_{s},\R^{+})>0,
$$
which is impossible, since the distribution of $X_{s}$ has a density with respect to Lebesgue measure.  From \eqref{E:Fnu}, it follows that $F^{(t_{n})}(\nu',\R^{+})=0$ for all compactly supported $\nu'$.  An identical argument, using \eqref{E:qom2} in place of \eqref{E:qom1}, shows the second part.
\end{proof}

\begin{Prop}
\label{P:collet}
Suppose the conditions (LP) and (GB) hold, and that $K:=\lim_{z\to\infty}\kappa(z)$ exists, and is not equal to $\ulm$.  Then the following are equivalent:
\begin{itemize}
\item[(a) ] For any sequence $t_{n}$ such that $F_{t_{n}}$ and $\omega_{t_{n}}$ converge, $F^{(t_{n})}(\nu,\R^{+})$ is constant in $\nu$.

\item[(b) ] If $\nu$ is any compactly supported initial distribution, then $\lim_{t\to\infty} F_{t}(\nu,\R^{+})$ exists.  For each $\nu$, the limit is either 0 or 1.

\item[(c) ] There is some $\nu\in\mc$ such that $\lim_{t\to\infty} F_{t}(\nu,\R^{+})$ exists, and is either 0 or 1.

\item[(d) ] For all $x\in\R^{+}$, the limit $\lim_{t\to\infty} \omega_{t}(x,\R^{+})$ exists. For each $x$, the limit is either $\psi_{K}$ or $\psi_{\ulm}$.
\end{itemize}
\end{Prop}

\begin{proof}
(a)$\Rightarrow$(b) Suppose that $F^{(t_{n})}(\nu,\R^{+})$ is constant.  By Proposition \ref{P:harmonic}, we know that $\omega^{(t_{n})}(x)$ is either $\ulm$-invariant or $K$-invariant, since $F^{(t_{n})}(x,\R^{+})$ and $1-F^{(t_{n})}(x,\R^{+})$ cannot both be 0.  But one of these must be 0, since $\omega^{(t_{n})}(x)$ cannot be both $\ulm$-invariant and $K$-invariant.  Thus $F^{(t_{n})}(\nu,\R^{+})$ is either 0 or 1.

Suppose, now, that there were a sequence $t_{n}$ such that $F^{(t_{n})}(\nu,\R^{+})=0$, and another sequence $t'_{n}$ such that $F^{(t'_{n})}(\nu,\R^{+})=1$.  Since $F_{t}(x,\cdot)$ is continuous in $t$, it would be possible to find an intermediate sequence $t''_{n}$ such that the limits $F^{(t''_{n})}$ and $\omega^{(t''_{n})}$ exist, and $F^{(t''_{n})}(x,\R^{+})\in (0,1)$.  This proves that $F(x,\R^{+})$ is sequence-independent.

\vspace{3mm}

(b)$\Rightarrow$(a) From the invariance properties of Proposition \ref{P:harmonic}
 it follows that $F(x,\R^{+})$ is constant in $x$.  Suppose now that $\ltf F_{t}(x,\R^{+})=0$, and let $\nu$ be any compactly supported probability.  Choose $\gep>0$ and $z>0$.  For each $x$, we may find $t(x)$ such that $F_{t}(x,[0,z])<\gep$ for all $t\ge t(x)$.  Thus
$$
P_{\nu}\ls X_{t}\le z \cond \tp>t \rs = \frac{\int \P_{x}\{X_{t} \le z \} d\nu(x)}{\int \P_{x}\{\tp>t\} d\nu(x)}
\le\nu\ls x\, :\, t(x)>t \rs + \gep,
$$
which converges to $\gep$ as $t\to\infty$.  Letting $\gep\to 0$ proves that $F_{t}(\nu,[0,z])=0$.  A similar proof holds if $F(x,\R^{+})\equiv 1$.

\vspace{3mm}

(b)$\Rightarrow$(c) This is trivial.

\vspace{3mm}

(c)$\Rightarrow$(a) Suppose first that $\ltf F_{t}(\nu,\R^{+})=0$.  Then $F^{(t_{n})}(\nu,\R^{+})=0$ for any sequence $t_{n}$. By Corollary \ref{C:solidarity}, it follows that $F^{(t_{n})}(\nu',\R^{+})=0$ for any $\nu'$.  An equivalent proof works if $\ltf F_{t}(\nu,\R^{+})=1$.

\vspace{3mm}

(b)$\Leftrightarrow$(d) This is a direct consequence of \eqref{E:omegaform}.
\end{proof}

\nc{\bP}{\bar{P}}
\subsection{The Dichotomy: Proving Theorem \ref{T:limsups} for $\mathbf{\ulm<\liminf_{z\to\infty} \kappa(z)}$}  \label{sec:bigK}
Assume $\ulm<\liminf_{z\to\infty} \kappa(z)$.  If $\izf \plm(y)dy=\infty$, then the process escapes to infinity from every starting distribution satisfying (IDC), by Proposition \ref{P:esc}.  Since (IDC) distributions are dense in $\mc$, by the continuity of $F$ we see that $F^{(t_{n})}(\nu,\R^{+})=0$ for any compactly supported $\nu$.

Assume now that $\izf \plm (y)dy$ is finite.  We need to show that
\begin{equation}  \label{E:needtoshow}
\lim_{z\to\infty}\ltf \frac{\P\{ X_{t}\ge z\}}{\P\{ X_{t}\le z\}}=0
\end{equation}
(We drop the compactly supported starting measure from the notation.)

    We first prove \eqref{E:needtoshow} with $t$ running only over 
    integers.
    For positive integers $n$ and positive real $z$, 
    we define 
    \begin{equation}  \label{E:fzn}
    f(z,n):=\frac{\Pr X_{n}>z\rs}{\Pr X_{n}\le z\rs}.
    \end{equation}    
    Choose $z_{0}$ such that $K:=\inf_{y\ge z_{0}} \kappa(y)>\ulm$.
    Observe that for $z>z_0$
    \begin{equation}
        \begin{split}\label{E:allterms}
        \Pr &X_{n+1}>z \rs \le \Pr X_{n}\le z \text{ and } X_{n+1}> 
        z \rs  \\
         & \hspace*{1cm} +\Pr X_{n}>z \text{ and } X_{n+1}>z \tand 
         X_{t}>z_{0} \; \fa t\in [n,n+1] \rs  \\
         & \hspace*{1.5cm} + \Pr X_{n}>z \text{ and } X_{n+1}>z \tand 
         \te t\in [n,n+1]\text{ s.t. }X_{t}\le z_{0} \rs .
        \end{split}
    \end{equation}
    Thus
    \begin{equation}
        \begin{split}\label{E:allterms2}
        f&(z,n+1)\le \Pr X_{n+1}>z \cond X_n\le z \rs  \\
        &+  f(z,n) \frac{\Pr X_n\le z \rs}{\Pr X_{n+1}\le z \rs}
        \Bigl[ 
         \Pr X_{n+1}>z \,\&\, 
         X_{t}>z_{0}, \; \fa t\in [n,n+1] \cond X_n>z\rs\\
         & \hspace*{1cm} + \Pr X_{n+1}>z \, \&\,
         \te t\in [n,n+1]\text{ s.t. }X_{t}\le z_{0} \cond X_n>z\rs 
         \Bigr].
        \end{split}
    \end{equation}
    By the strong Markov property (taking the first time the process hits 
    the level $z_0$ after time $n$), the last term in brackets is 
    bounded above by
    \begin{equation}
        \mathbb{P}_{z_{0}} \ls \te t\in [0,1] \tst 
        X_{t}>z \rs.
        \label{eq:lastterm}
    \end{equation}
    The middle term in the brackets in \eqref{E:allterms2} is bounded by $e^{-K}$, since the process is being killed at rate at least $K$ while it is 
    above $z_{0}$.
        
    Substituting these bounds into \eqref{E:allterms2}, we get
     \begin{equation}\label{E:allterms3}
        \begin{split}
        f(z&,n+1)\\
        & \le\frac{\Pr X_n\le z \rs}{\Pr X_{n+1}\le z \rs}
        \Bigl[\Pr X_{n+1}> z \cond X_{n}\le 
        z \rs\\
        &\hspace*{2cm}+ \lp e^{-K} + \mathbb{P}_{z_{0}} \ls \te t\in [0,1] \tst 
        X_{t}>z \rs \rp f(z,n) \Bigr] .
        \end{split}
    \end{equation}
    We write this in the form
    $$
    f(z,n+1)\le c(z,n) f(z,n)+\gm(z,n),
    $$
    where
    \begin{equation}
        c(z,n):=\lp \mathbb{P}_{z_{0}} \ls \te 0\le t\le 1 \, : \,
        X_{t}>z \rs +e^{-K}\rp\Bigl( \frac{\Pr X_n\le z \rs}{\Pr X_{n+1}\le z \rs}
        \Bigr),
        \label{eq:cnz}
    \end{equation}
    and
    \begin{equation}
        \gm(z,n):=\Pr X_{n+1}> z \cond X_{n}\le
        z \rs \Bigl( \frac{\Pr X_n\le z \rs}{\Pr X_{n+1}\le z \rs}\Bigr)^{-1}.
        \label{eq:gmnz}
    \end{equation}
    
    From Lemma \ref{L:ratiowork}, we know that
    $$
    \lnf \frac{\Pr X_n\le z \rs}{\Pr X_{n+1}\le z \rs} = e^{\ulm}.
    $$
Substituting into \eqref{eq:cnz} yields,
    \begin{equation*}
        \begin{split}
            \lzn c(z,n)&\le e^{(\ulm-K)}+e^{\ulm} \lzn \int_{0}^{z}
            \mathbb{P}_{z_0}\ls \te 0\le t\le 1 \tst 
        X_{t}>z \rs\\
            &=e^{(\ulm-K)}.
        \end{split}
    \end{equation*}
    (The second term on the right is 0 because $\infty$ is an inaccessible 
    boundary.)  
    
    By quasistationary convergence on compacta,
    $$
    \lzn \gm(z,n)\le \limsup_{z\to\infty}\int_{0}^{z}
            \mathbb{P}_{y}\ls X_{1}>z \rs d\mu_{z}(y),
    $$
    where $\mu_{z}$ is the distribution on $[0,z]$ with density proportional to $\phi_{\ulm}$.  Since $\infty$ is inaccessible, it follows that $\lim_{z\to\infty}\mathbb{P}_{y}\ls X_{1}>z \rs=0$ for each fixed $y$, so that $\lzn \gm(z,n)=0$.
    Lemma \ref{L:twolims} implies then that
    $$
    \lzn f(z,n)=0.
    $$

Now we need to extend the proof of \eqref{E:needtoshow} to real 
times $t$.  We have
\begin{equation}  \label{E:needtoshow2}
\begin{split}
\lim_{z\to\infty} &\limsup_{t\to\infty} \frac{\Pr X_{t}>z\rs}{\Pr X_{t}\le z\rs} \\
&\le
\lim_{z\to\infty} \limsup_{n\to\infty} \frac{\Pr \te t\in [n,n+1) \text{ s.t. }
X_{t}>z\rs}{\P\{ X_{n}\le z\}} \sup_{t\in [n,n+1)}\frac{\P\{ X_{n} \le z\}}{\P\{ X_{t} \le z\}}.
\end{split}
\end{equation}
The numerator of the first term on the right is bounded above by
$$
\Pr X_{n}>z \rs + \Pr X_{n}\le z \rs \cdot \Pr \te  t\in [n,n+1) \text{ s.t. }
X_{t}>z \cond X_{n}\le z \rs.
$$
By quasistationary convergence on compacta, the limit as $n\to\infty$ of the right-hand term is
$$
\left(\int_{0}^{z}\phi_{\ulm}(x)dx\right)^{-1}
\int_{0}^{z} \mathbb{P}_{y}\ls \te t\in [0,1) \text{ s.t. }X_{t}>z \rs 
\phi_{\ulm}(y)dy.
$$
Since $\infty$ is inaccessible, this goes to 0 as $z\to\infty$.

The expression in \eqref{E:needtoshow2} is then bounded by
\begin{multline*}
\lzn \sup_{t\in [n,n+1]} \frac{\P\{ X_{n} \le z\}}{\P\{ X_{t} \le z\}}
\frac{\Pr X_{n}>z\rs}{\Pr X_{n}\le z \rs},
\end{multline*}
which we have already shown to be 0.

\nc{\bk}{\bar{K}}
\subsection{The Dichotomy: Proving Theorem \ref{T:limsups} for $\mathbf{\ulm>\limsup_{z\to\infty} \kappa(z)}$}  \label{sec:smallK}
Heuristically, the proof is very simple, though the notation may be confusing.  Consider some large $z$, large enough that the killing rate is always smaller than $\ulm$ above $z$; and choose some $z'>z$, large enough that the probability of moving from above $z'$ to below $z$ in one unit of time is small.  Eventually, mass below $z$ must be declining at rate $\ulm$.  The mass above $z'$ is being killed at a strictly lower rate.  It could be, then, that the mass above $z'$ blows up relative to the mass down below, in which case we have escape to infinity.  Otherwise, the mass above $z'$ must be sliding down.  Can this happen?  The obstruction is the narrow pipe between $z$ and $z'$.  In one unit of time, only a small fraction of the mass above $z'$ can make its way down below $z$.  Thus, if mass is escaping downward, it must be spending some time in $[z,z']$.  Convergence to a probability distribution on compacta implies, though, that eventually the mass on $[z,z']$ --- taking up a significant fraction of mass from $[z',\infty)$ --- will be a small fraction of the mass on $[0,z]$.  Putting all of this together, we conclude that the mass above $z'$ must be much smaller than that below $z$.

We now translate the heuristic into formal mathematics.
Fix $z$ large enough that
$$
\sup_{y\ge z} \kappa(y)\le (K+\ulm)/2=:\bk,
$$
and define 
$$
\gep:=\frac{\int_{z}^{\infty}\plm(y)dy}{\int_{0}^{z}\plm(y)dy}.
$$
Associate to $z$ some $z'(z)>z$ (assumed increasing in $z$) such that
$$
\P_{z'}\ls \te \, t\in[0,1]\, :\, X_{t}\le z\rs<\gep .
$$
(This is possible, since we have assumed that $\infty$ is a natural boundary.)  We define
\begin{align*}
a_{n}:=\P\ls X_{n}\le z \rs,\\
b_{n}:=\P\ls z\le X_{n}\le z' \rs,\\
c_{n}:=\P\ls X_{n}\ge z' \rs.
\end{align*}
We have the following relations:
\begin{align*}
c_{n+1}&\ge e^{-\bk}c_{n}\Bigl(1-\P_{z'}\ls \te \, t\in[0,1]\, :\, X_{t}\le z\rs\Bigr)-b_{n+1},\\
\limsup_{n\to\infty} \frac{b_{n}}{a_{n}}&\le \gep,\\
\lnf \frac{a_{n+1}}{a_{n}} &=e^{-\ulm}.
\end{align*}
Thus,
$$
\frac{c_{n+1}}{a_{n+1}}\ge \frac{a_{n}}{a_{n+1}} e^{-\bk}(1-\gep) \frac{c_{n}}{a_{n}} - \frac{b_{n+1}}{a_{n+1}}.
$$
Set $R=e^{\ulm-\bk}(1-\gep)^{2}$.  We may find $N$ such that for all $n\ge N$,
$$
\frac{c_{n+1}}{a_{n+1}}\ge R \frac{c_{n}}{a_{n}} - 2\gep.
$$

Suppose $z$ large enough so that $R>1$.  Then either
$$
\lim_{n\to\infty}\frac{c_{n}}{a_{n}} =\infty\quad \text{or }
\limsup_{n\to\infty}\frac{c_{n}}{a_{n}}\le \frac{2\gep}{R-1}.
$$
These translate into
$$
\lim_{n\to\infty} f(z,n) =\infty \quad \text{or}\quad
\limsup_{n\to\infty} f(z,n)\le  \frac{2\gep}{R-1} +\gep
$$
for $z$ sufficiently large, where $f(z,n)$ is as defined in \eqref{E:fzn}.
Note that $f(z,n)$ is monotone decreasing as a function of $z$, so that if the second limit holds for any $z$, it holds for all higher $z$. Sending $z$ to $\infty$ makes $\gep\to 0$, so that either
\begin{equation} \label{E:eitheror}
\lim_{n\to\infty} f(z,n) =\infty \text{ for all }z\quad \text{or}\quad
\lzf\limsup_{n\to\infty} f(z,n)=0.
\end{equation}

The same argument as in section \ref{sec:bigK} shows that the second condition in \eqref{E:eitheror} is equivalent to
$$
\lim_{z\to\infty} \limsup_{t\to\infty} f(z,t)=0.
$$
Similarly, if $\lnf f(z,n)=\infty$ for all $z$, we have for any $z'>z$,
\begin{align*}
\liminf_{t\to\infty} &f(z,t) \ge \liminf_{n\to\infty}\frac{ \P\{ X_{n}>z' \} - \P\{ X_{n}>z' \; \& \; \te t\in (n,n+1) \, :\, X_{t}\le z \}}{ \P\{ X_{n}\le z' \} + \P\{ X_{n}>z' \, \& \, \te t\in (n,n+1) \; :\; X_{t}\le z \}}\\
&\ge  \liminf_{n\to\infty} \frac{f(z',n)(1- \P \{ \te t\in (n,n+1) \, :\, X_{t}\le z \cond X_{n}>z'\})}{1+ f(z',n)\P \{ \te t\in (n,n+1) \, :\, X_{t}\le z \cond X_{n}>z'\}}\\
&=\liminf_{n\to\infty} \P \{ \te t\in (n,n+1) \, :\, X_{t}\le z \cond X_{n}>z'\}^{-1} -1.
\end{align*}
Since $\infty$ is a natural boundary, this diverges to $\infty$ as $z'\to\infty$.  Thus, we have shown that
\begin{equation}  \label{E:eitheror2}
\lim_{t\to\infty} f(z,t)=\infty\text{ for all } z \quad\text{or}\quad
\lim_{z\to\infty} \limsup_{t\to\infty} f(z,t)=0.
\end{equation}

If we were assuming (GB), we could apply Proposition \ref{P:collet} here.  Instead, we note that $\nu\mapsto F^{(t_{n})}(\nu,\R^{+})$ is continuous.  Since we have shown that it takes on only the values 0 and 1, it must be constant.  The case $F^{(t_{n})}(\nu,\R^{+})\equiv 0$ is the same as escape to infinity; the case $F^{(t_{n})}(\nu,\R^{+})\equiv 1$ is (by Lemma \ref{L:whichdensity}) the same as convergence to the quasistationary distribution $\plm$.

\subsection{Proof of Theorem \ref{T:limits}}  \label{sec:noescape}
If $\ulm=K$, then Lemma \ref{L:omegat} together with \eqref{E:omegaform} implies that every sequence $t'_{n}\to\infty$ has a subsequence $t_{n}$ such that $\lnf \omega_{t_{n}}(x)=\psi_{\ulm}(x)$.  Thus, $\psi_{\ulm}$ is the limit.  The same argument applies to $a_{t_{n}}(x,r)$, proving the claim about the killing rate.

Suppose now that $\ulm<K$.  Since $\ulm$ is the bottom of the spectrum of $\mL$, it follows (as we have made explicit in Lemma \ref{L:bottom}) that there is some $x>0$ for which $\psi_{K}(x)=0$.  Since $\omega$ is strictly positive (except perhaps at 0), it follows from \eqref{E:qom2} that $F^{(t_{n})}(x,\R^{+})=1$ for any sequence $(t_{n})$ for which the limits exist and are continuous.  By Proposition \ref{P:collet},  $\ltf F_{t}(\nu,\R^{+})=1$ for every compactly supported $\nu$.  The results of Theorem \ref{T:limits} follow then from Lemma \ref{L:whichdensity} and Lemma \ref{L:arelate}.

Suppose, finally, that $\ulm>K$.  If $X_{t}$ escapes to infinity from any compactly supported initial distribution,   it follows, by Proposition \ref{P:collet}, that $X_{t}$ escapes to infinity from every compactly supported initial distribution.  Suppose, then, that $X_{t}$ does not escape to infinity from any compactly supported initial distribution.  Then $\izf \plm(y)dy<\infty$, so we may apply Theorem \ref{T:limsups} to see that  $\ltf F_{t}(\nu,\R^{+})=1$ if $\nu$ satisfies (IDC).  By Proposition \ref{P:collet}, this holds then for all compactly supported $\nu$.  The result then follows from Lemma \ref{L:whichdensity}.

\section{An example}  \label{sec:example}
One Markov mortality model that has received considerable attention is 
a discrete-state-space process that was first defined by H. Le Bras 
\cite{hLB76}.  (A somewhat similar process appeared in the mathematical literature, presumably independently, in \cite{jC78}.)  In this process the states are positive integers, and the motion is only increasing by single steps.  The rate of transition 
to the next higher state, and the mortality rate, are both directly proportional to the current state.

We define a continuous-state-space analogue to this process on $[1,\infty)$ by the stochastic differential equation
\begin{equation}
    \label{E:lebrasdif}
    dX_{t}=\gs X_{t} dW_{t} + b X_{t} dt,
\end{equation}
where $\gs$ is a positive constant and b is a constant larger than
$\gs^{2}/2$.  The process starts at $X_{0}=1$, is killed at the rate
$k X_{t}$, and is reflected when it hits 1.  By Ito's formula we see
that $X_{t}$ has the same dynamics in $(1,\infty)$ as the geometric 
Brownian motion
$$
X_{t}= \exp\left\{ \gs (W_{t} + \tilde{b} t) \right\},
$$
where $\tilde{b}=\frac{b}{\gs}-\frac{\gs}{2}$.  Equivalently, then, we could consider the 
Brownian motion with drift:
$$
Y_{t}=W_{t}+ \tilde{b} t,
$$
killed at a rate $ke^{\gs y}$ and reflected at 0.  This will have 
the same mortality distribution as $X_{t}$.  We have $B(x)=\tilde{b}x$, where 
$\tilde{b}$ is positive.  The left boundary 0 is regular, while the right 
boundary $\infty$ is natural.  Since $\kappa(x)\to\infty$ as 
$x\to\infty$, we can apply Theorem \ref{T:limsups}.  This guarantees the convergence to a quasistationary distribution from any initial distribution with a compactly supported density, as long as $\izf \plm (y)dy<\infty$.  Condition (GB) is satisfied, as we discuss in the comments preceding Lemma \ref{L:omegabounded}, since the drift is constant and $\kappa$ is increasing.  This allows us to apply Theorem \ref{T:limits}, to conclude that $\izf \plm (y)dy<\infty$.  (This would follow from our computation in any case.)

Eigenfunctions $\phi_{\gl}$ for the generator of $Y_{t}$ are given (up to a scaling constant) by
$$
\phi_{\gl}(x)=\gs e^{\gs y} \xi_{\gl}(\gs e^{\gs y}),
$$
where $\xi_{\gl}:[1,\infty)\to\R$ satisfies
\begin{equation}
    \oh \gs^{2} \lp x^{2} \xi_{\gl}(x) \rp'' - b\lp x \xi_{\gl}(x)\rp'
       -kx \xi_{\gl}(x) = -\gl \xi_{\gl}(x) , \label{E:fklb}
\end{equation}
with the boundary condition
\begin{equation}  \label{E:bdryphi}
    \xi_{\gl}'(1) = \left( \frac{b}{\gs^{2}} - \frac{3}{2} \right) \xi_{\gl}(1).
    \end{equation}
In a sense, $\xi_{\gl}$ is the function of interest, since this is the density of the quasistationary distribution in the original coordinates.  Formally, though, we need to refer to $\phi_{\gl}$, since our theorems are stated in this normalization.  In particular, it is $\plm$ that needs to be integrable.)

If we define a new function $\psi$ by
$$
\xi_{\gl}(x)=x^{\frac{b}{\gs^{2}}-\frac{3}{2}}
          \psi\left(\sqrt{\frac{8kx}{\gs^{2}}}\right),
$$
we get
$$
z^{2}\psi''(z)+z\psi'(z)-(\nu^{2}+z^{2})\psi(z) =0,
$$
where
\begin{equation}  \label{E:nu}
\nu= \sqrt{\left( \frac{2b}{\gs^{2}} - 1\right)^{2} - 
\frac{8\gl}{\gs^{2}}}.
\end{equation}
(Note: In order to find the transformation, we begin by writing 
$\xi_{\gl}(x)=x^{r}\psi(\sqrt{ax})$, and solve for $a$ and $r$.)
This is the equation satisfied by the modified 
Bessel functions.  Thus the general solution for $\phi$ is
$$
\xi_{\gl}(x)=x^{\frac{b}{\gs^{2}}-\frac{3}{2}}
\left[C_{1}(\gl)I_{\nu}\left(\sqrt{\frac{8kx}{\gs^{2}}}\right) \, + \, 
C_2(\gl) K_{\nu}\left(\sqrt{\frac{8kx}{\gs^{2}}}\right) \right],
$$
where $\nu$ is 
given as a function of $\gl$ by \eqref{E:nu}.  On any compact set of $z$ 
and $\gl$, $\xi_\gl(z)$ is uniformly continuous in $z$ and $\gl$ (by Theorem
1.7.5 of \cite{CL55}), which implies that $C_1$ and $C_2$ are continuous.

\nc{\ula}{\underline{\gl}}

We need to compute $\ulm$.  By Lemma \ref{L:bottom}, we know that $\xi_{\ulm}$ is positive, but $\xi_{\gl}$ changes sign for all $\gl>\ulm$   As $z\to \infty$, $I_{\nu}(z)$ is increasing like 
$e^{z}/\sqrt{z}$ and $K_{\nu}$ decreasing like $e^{-z}/\sqrt{z}$.
Togther with the continuity of the coefficients, this implies that $C_1(\ula)=0$. The top eigenfunction that we seek is a multiple of $I_\nu$.

The transformation of $\xi_{\gl}$ to $\psi$ shifts the boundary from 1 to 
$$x_{0}=\sqrt{8k}/\gs.
$$  
We have
\begin{equation}\label{E:otherbound}
    \frac{\phi'(x)}{\phi(x)} = \left( 
    \frac{b}{\gs^{2}}-\frac{3}{2} \right)x^{-1}  +
    x^{-1}\frac{\sqrt{2kx}}{\gs}\frac{K'_{\nu}\left(\sqrt{8kx}/\gs\right)}
    {K_{\nu}\left(\sqrt{8kx}/\gs\right)}.
\end{equation}
Thus, the boundary condition turns into $K'_{\nu}(x_{0})=0$.

There cannot be a solution to 
\eqref{E:otherbound} for $\nu$ real and positive, since then 
$K_{\nu}(x_{0})$ is positive, while the derivative is negative.  By 
the integral representation (9.6.24) of \cite{AS65},
\begin{align*}
    K_{iy}(x) & = \izf e^{-x\cosh t} \cos(y t)dt, \\
    K'_{iy}(x) & = -\izf e^{-x\cosh t} \cosh(t)\cos(y t)dt. 
\end{align*}
Observe, first, that $K_{iy}(x)$ has infinitely many zeros close to 
$x=0$, but also a final zero; the location of this final zero is 
monotonically increasing in $y$, and unbounded.  

The same is true for 
the derivative $K'_{iy}(x)$, and the final zero of $K'$ is 
larger than the final zero of $K$.  The zeros of $K$ and $K'$ are 
interleaved.  For $y=0$, 
$K_{iy}(x_{0})$ is positive, while $K'_{iy}(x_{0})$ is negative.
For $y$ close to 0 this ratio remains negative.  As $y$ increases, 
eventually the final zero of $K'_{iy}$ reaches $x_{0}$, at $y=\tilde{y}(x_{0})$.  (Remember that $x_{0}$ is itself a function of the process parameters $\gs$ and $k$.)  The limiting rate of killing is then
$$
\ulm=\gl(i\tilde{y})=\frac{\gs^{2}}{8}\left[\left(\frac{2b}{\gs^{2}}-1 \right)^{2} 
\, + \, \tilde{y}^{2}\right].
$$
The quasistationary distribution is a constant multiple of
$$
\xi_{\ulm}=x^{\frac{b}{\gs^{2}}-\frac{3}{2}}K_{i \tilde{y}}
  \left(\sqrt{\frac{8kx}{\gs^{2}}}\right).
$$
This function behaves asymptotically as
$$
x^{\frac{b}{\gs^{2}}-2} e^{-\sqrt{8kx/\gs^{2}}}
$$
as $x\to\infty$, according to formula (9.7.2) of \cite{AS65}.
We see that $\plm(x)=\gs e^{\gs x}\xi_{\ulm}(e^{\gs x})$ is integrable, confirming the conclusion of Theorem \ref{T:limits}.

\section{Technical lemmas}  \label{sec:tech}
\subsection{Facts about generators and resolvents}  \label{sec:gr}
\begin{Lem}  \label{L:gen}
The adjoint semigroup for the diffusion with killing is generated by the operator with domain $\mD=\{ f \in L^{1}\cap \mC\, :\, \mL^{*} f \in L^{1}\},$ acting as $\mL^{*}$.
\end{Lem}

\begin{proof}
The Feynman-Kac construction already defines a strongly continuous contraction semigroup and an adjoint strongly continuous contraction semigroup. Initially, the adjoint semigroup acts on $\eu{BV}([0,r])$.  We observe, first, that the image of any measure under the resolvent has no atoms at the endpoints, since we have assumed instantaneous reflection (in the case of regular boundaries) or instantaneous entrance (in the case of entrance boundaries).  Now, let $R_{\ga}^{k}$ be the resolvent for the semigroup with internal killing, and $R_{\ga}^{k*}$ the adjoint resolvent with internal killing.  For $\nu$ any finite measure on $\R^{+}$ and $A\subset \R^{+}$,
\begin{align*}
R_{\ga}^{*}\nu(A) &= \E_{\nu} \izf e^{-\ga t} \indic\{X_{t}\in A\} dt,\text{ and}\\
 R_{\ga}^{k*}\nu(A) &= \E_{\nu} \izf e^{-\ga t} e^{-\int_{0}^{t}\kappa(X_{s})}\indic\{X_{t}\in A\} dt.
\end{align*}
We see then that $R_{\ga}^*\nu $ and $R_{\ga}^{k*}\nu $ are equivalent measures.  Since we know from Theorem 13.3 of \cite{wF52} that $R_{\ga}^{*}\nu$ has a density, so must $R_{\ga}^{k*}\nu$ have a density.  Both are bounded in total variation by $1/\ga$.  This tells us, in particular, that all of $L^{1}$ is mapped into $L^{1}$ by the resolvent.  Since the range of the resolvent is dense in the range of the semigroup, this means that the restriction of the adjoint semigroup action to $L^{1}$ is a strongly continuous contraction semigroup on $L^{1}$.

Call the generator of this adjoint semigroup with internal killing $\Omega^{*}$. Formally, the adjoint generator must act as the formal adjoint of the generator, namely, as $\mL^{*}$, on functions in the domain which are twice differentiable.  It remains to show that $\mD$ is a core for $\Omega^{*}$.  By Proposition 1.3.1 of \cite{EK86}, it suffices to show that $\mD$ is a dense subset of the range of the resolvent, and that the resolvent maps $\mD$ into itself.

We first need to show that $\mD$ is actually contained in the range of the resolvent.  That is, for every $f\in\mD$ we must find $g\in L^{1}$ such that $R^{k*}_{\ga} g=f$. Equation 6 of section 4.11 of \cite{IM65} represents the $\ga$-Green's function (which is the resolvent density) for the killed diffusion by
\begin{equation}  \label{E:green}
\G_{\ga}(y,x)= \G^*_\ga(x,y)=\begin{cases}
    cg_{1}(x)g_{2}(y)e^{B(x)} &\tif x\le y,\\
    cg_{1}(y)g_{2}(x)e^{B(x)} &\tif y\le x,
    \end{cases}
\end{equation}
where $g_{1}$ and $g_{2}$ are increasing and decreasing
eigenfunctions, respectively, for $\mL$, and $c$ is a constant. Define $g=\ga f - \mL^{*}f$ and $h=R_{\ga}^{k*}g$.  By the definition of $\mD$, $f$ is in
$L^{1}$ and twice continuously differentiable.  Thus $g$ is
continuous, and it follows, by Lemma \ref{L:grapplies}, that $h$ is continuously twice differentiable, so that $\Omega^{*}$ acts as the differential operator $\mL^{*}$ on h.  This implies that both $f$ and $h$ are solutions to the differential equation
$$
\ga u - \mL^{*} u = g,
$$
so that
$$
\mL^{*}[f-h]=\ga[f-h].
$$

We show first that $h$ satisfies the boundary conditions \eqref{E:fellerbound} at $0$.
The boundary condition for the image of $R_{\ga}^{k}$ is given on page 131 of \cite{IM65} as
\begin{equation}  \label{E:fellerbound2}
(1-p_{0}) \lim_{x\to 0}\phi(x) = \oh p_{0}\lim_{x\to 0} e^{B(x)} \phi'(x).
\end{equation}
This is equivalent to the equation for the Green's function
\begin{equation}  \label{E:greenbound}
\fa y,\quad (1-p_{0}) \lim_{x\to 0}G_{\ga}(x,y) = \oh p_{0}\lim_{x\to 0} e^{B(x)} \frac{\partial G_{\ga}(x,y)}{\partial x}.
\end{equation}
Now, from the representation of the Green's function as \eqref{E:green}, we see that $G_{\ga}(x,y)=e^{B(y)-B(x)}G_{\ga}^{*}(x,y)$.  This means that
$$
\lim_{x\to 0} e^{B(x)} \frac{\partial G_{\ga}(x,y)}{\partial x}=e^{B(y)}\lim_{x\to 0}\left[ \oh \frac{\partial G_{\ga}^{*}(x,y)}{\partial x} - b(x) G_{\ga}(x,y) \right].
$$
This translates into the condition on $h=R_{\ga}^{k*}g$ of
$$
(1-p_{0}) \lim_{x\to 0}e^{-B(x)}h'(x) = p_{0} \lim_{x\to 0} \lb h'(x) - 2b(x) h'(x) \rb,
$$
which is exactly the condition in \eqref{E:fellerbound}.  Similarly, if $r$ is regular, $h$ satisfies the condition \eqref{E:fellerboundr}, and if $r$ is an entrance boundary $h$ satisfies \eqref{E:fellerboundentrance}.  Thus $f-h$ is an $L^1$ solution to $\mL^*u=\ga u$, satisfying these boundary conditions.  We show in Lemma \ref{L:eigenfunction} (copying \cite{wF52}) that this implies that $f-h=0$.

We now show that $R^{k*}_{\ga}$ maps $\mD$ to itself.  We have just shown that the adjoint resolvent takes all continuous $L^{1}$ functions to continuously twice-differentiable functions which satisfy the appropriate boundary conditions.  If $f=R^{k*}_{\ga}g$, where $g\in \mD$, then $f\in\mC$, and
$$
\mL^{*}f=\ga f - g \in L^{1}.
$$
Finally, if $f$ is in the domain of $\Omega^{*}$, we can represent $f$ as $R^{k*}_{\ga}g$ for some $g\in L^{1}$, then we can represent $g$ as a limit of continuous functions $g_{i}$.  Then $f_{i}=R^{k*}_{\ga}g_{i}$ is in $\mD$, and satisfies $\|f_{i}-f\|\le \ga^{-1} \|g_{i}-g\|$, implying that $\lim_{i\to\infty} f_{i}= f$.  As we have already shown that the functions $g_{i}$ are in $\mD$, this shows that $\mD$ is dense in the domain, completing the proof.
\end{proof}

\nc{\god}{\Gamma^{\Delta}_{1}}
\nc{\gD}{\Gamma^{\Delta}}

\begin{Lem}
\label{L:grapplies}
Letting $G^{*}_{\ga}(x,y)$ be the adjoint Green's function given in \eqref{E:green}, and $g$ any continuous $L^{1}$ function, then
$$
h(x):=\int_{0}^{r} G^{*}_{\ga}(x,y)g(y)dy
$$
is continuously twice differentiable.
\end{Lem}

\begin{proof}
The technical center of this result is a merely tedious calculus fact, based on the differentiability of $G^{*}$ away from the diagonal.  We isolate the details as Lemma \ref{L:green}.  It remains only to show that the conditions of that lemma are satisfied by $G^{*}$ for $n=2$.

Since we assumed that $b$ is continuously differentiable, and that $\kappa$ is continuous, it follows that $B$ is twice continuously differentiable, as are the eigenfunctions $g_{1}$ and $g_{2}$ of $\gL^{*}$.
We have $G^{*}(x,y)=g_{1}(x\wedge y)g_{2}(x\vee y)e^{B(x)}$.  We have
$$
G^{*}_{2}(x,y)=\begin{cases}
c(g_{2}(x)e^{B(x)})'' g_{1}(y) &\tif y<x;\\
c(g_{1}(x)e^{B(x)})'' g_{2}(y) &\tif y>x.
\end{cases}
$$
Since $g_{1}$ is increasing and $g_{2}$ decreasing, the factor which depends on $y$ is bounded.  This implies that $G^{*}_{2}(x,y)$ and $G^{*}_{1}(x,y)$ are bounded in $y$, and that their continuity and convergence are uniform in $y$.  As for $G^{*\Delta}_{1}$, we have
$$
G^{*\Delta}_{1}(x)=\lp g_{1}(x) e^{B(x)}\rp' g_{2}(y)- \lp g_{2}(x) e^{B(x)}\rp' g_{1}(y),
$$
which is continuous.  Since $G^{*}$ is continuous, it follows that $G^{*\Delta}_{0}$ is identically 0.
\end{proof}

\begin{Lem}
\label{L:green}
Let $\Gamma(x,y)$ be a real-valued function on $[0,r]\times [0,r]$, which is $n$-times differentiable as a function of $x$ for every interior point $y\ne x$, where $n$ is a positive integer.  Let $g\in L^{1}(0,r)$ be continuous, and define 
$$
h(x):=\int_{0}^{r} \Gamma(x,y)g(y)dy.
$$  
Define $\Gamma_{k}(x,y)=\partial^{k} \Gamma/\partial x^{k}(x,y)$ for $k\le n$, and let
$$
\Gamma^{\Delta}_{k}(y):=\lim_{x\downarrow y}\Gamma_{k}(x,y)-\lim_{x\uparrow y}\Gamma_{k}(x,y).
$$

Suppose that
\begin{itemize}
\item $\Gamma_{k}(x,y)$ is bounded as a function of $y$, for $k\le n$.  That is,
\begin{equation}  \label{E:nbounded}
\sup_{y\ne x}\lv\Gamma_{k}(x,y)\rv<\infty
\end{equation}

\item $\Gamma_{k}(x,y)$ is continuous in $x$, uniformly in $y\ne x$, for $0\le k\le n$.  That is,
\begin{equation}  \label{E:unifcont}
\begin{split}
\lim_{\gep\to 0^{+}}& \sup_{y\in (0,x)\cup (x+\gep,r)} \lv \Gamma_{k} (x+\gep,y)-\gep \Gamma_{k}(x,y) \rv =0, \tand\\
        \lim_{\gep\to 0^{+}}& \sup_{y\in (0,x-\gep)\cup (x,r)} \lv \Gamma_{k}(x-\gep,y)-\gep \Gamma_{k}(x,y) \rv =0.
\end{split}
\end{equation}
        \item The convergence to the $k$-th derivative is uniform in $y\ne x$, for $1\le k\le n$.  That is,
\begin{equation} \label{E:unifconv}
\begin{split}  
\lim_{\gep\to 0^{+}}& \sup_{y\in (0,x)\cup (x+\gep,r)} \gep^{-1}\lv \Gamma_{k-1} (x+\gep,y)-\Gamma_{k-1}(x,y) - \gep \Gamma_{n}(x,y) \rv =0, \tand\\
        \lim_{\gep\to 0^{+}}& \sup_{y\in (0,x-\gep)\cup (x,r)} \gep^{-1}\lv \Gamma_{k-1}(x-\gep,y)-\Gamma_{k-1}(x,y) - \gep \Gamma_{k}(x,y) \rv =0.
\end{split}
        \end{equation}
\item $\gD_{n-1}$ is continuous, and $\gD_{k}$ is zero for $k\le n-2$.
\end{itemize}

Then $h$ is $n$-times continuously differentiable on $(0,r)$.  The $n$-th derivative is given by
\begin{equation}  \label{E:deriv}
h_{n}(x):=f(x)\gD_{n-1}(x)+\int_{0}^{r} \Gamma_{n}(x,y) f(y)dy.
\end{equation}
\end{Lem}

\nc{\pg}{\frac{\partial \Gamma(x,y)}{\partial x}}
\begin{proof}
The proof is by induction on $n$.  First take $n=1$.  We have assumed that $\gD_{0}$ is continuous, and $\Gamma_{1}$ is continuous in $x$, uniformly in $y$.  Thus $h_{1}$ exists, and is continuous.
We need to show that for each $x\in (0,r)$,
\begin{align*}
0&=\lim_{\gep\to 0} \gep^{-1}\lb h(x+\gep)-h(x)-\gep h_1(x) \rb \\
&=\lim_{\gep\to 0} \gep^{-1} \left[ \int_{0}^{r}g(y)\left\{ \Gamma(x+\gep,y)-\Gamma(x,y) - \gep \Gamma_{1}(x,y) \right\} dy \; - \; \gep g(x)\gD_{0}(x) \right].
\end{align*}
We consider limit from the right; the limit from the left is essentially the same.  For $y$ outside of $[x,x+\gep]$, the quantity in braces is $o(\gep)$, uniformly in $y$.  This leaves us with the task of showing that
$$
0=\lim_{\gep\downarrow 0} \gep^{-1} \left[ \int_{x}^{x+\gep}g(y)\left\{ \Gamma(x+\gep,y)-\Gamma(x,y) - \gep \Gamma_{1}(x,y) \right\} dy  \; - \; \gep g(x)\gD_{0}(x) \right].
$$
We can rewrite 
\begin{align*}
\Gamma(x+\gep,y)-\Gamma(x,y)&=\int_{x}^{y}\Gamma_{1}(z,y) dz \; + \; \gD_{0}(y)  \; +\; \int_{y}^{x+\gep}\Gamma_{1}(z,y) dz \\
&= \int_{x}^{x+\gep}\Gamma_{1}(z,y) dz \; + \; \gD_{0}(y) ,
\end{align*}
transforming the limit into
 $$
 0=\lim_{\gep\downarrow 0} \gep^{-1} \int_{x}^{x+\gep}\left(g(y)\left[\int_{x}^{x+\gep} \left\{\Gamma_{1}(z,y)-\Gamma_{1}(x,y)\right\}dz\; +\; \gD_0(y)\right]\; -\; g(x)\gD_0(x) \right)dy.
 $$
It will suffice, then, to show that
$$
0=\lim_{\gep\downarrow 0}\sup_{y\in [x,x+\gep]}g(y)\left[\int_{x}^{x+\gep} \left\{\Gamma_{1}(z,y)-\Gamma_{1}(x,y)\right\}dz\; +\; \gD_0(y)\right]\; -\; g(x)\gD_0(x).
$$
By the assumption \eqref{E:nbounded}, boundedness of $\Gamma_{1}$, the integral is $O(\gep)$.  Since $g$ and $\gD_0$ are continuous, the remaining difference goes to 0 with $\gep$. This completes the proof for $n=1$.

Now, take $n\ge 2$, assume the lemma true for $n-1$, and that the conditions are satisfied for $n$.  If we substitute $\Gamma_{1}$ for $\Gamma$, the conditions are satisfied for $n-1$.  Furthermore, since $\gD_{0}\equiv 0$, we have
$$
h'(x)=h_{1}(x)=\int_{0}^{r}\Gamma_{1}(x,y)g(y) dy.
$$
It follows from the induction hypothesis that $h'(x)$ is $(n-1)$-times continuously differentiable, which completes the induction.
\end{proof}

The following lemma is copied from the results of sections 9 through 11 of \cite{wF52}, adapted for the case of nonzero internal killing.  With only a few small exceptions, the arguments are identical, though we have specialized them somewhat for the problem at hand.

\begin{Lem}
        \label{L:eigenfunction}
        For $\ga>0$, if $\phi\in L^1(0,r)$ is a solution to the differential equation 
        $\mL^* \phi=\ga \phi$ satisfying the boundary conditions \eqref{E:fellerbound}, and \eqref{E:fellerboundr} if $r$ is regular, then $\phi\equiv 0$.
\end{Lem}

\begin{proof}
Let $u$ be any solution to the forward eigenvalue problem $\mL u=\ga u$ on $(0,r)$.  We note first that $u$ cannot have any positive interior maxima.  If $x$ is a local maximum, then $0\ge u''(x)=2(\kappa(x)+\ga)u(x)$, which implies that $u(x)\le 0$. Similarly, there is no negative minimum. It follows that $u$ has limits at both endpoints.  It also follows that two solutions $u_{1}$ and $u_{2}$ with $u_{1}(x_0)=u_{2}(x_{0})=1$ do not meet at any other point.  (Recall that $x_{0}$ is an arbitrary point in $(0,r)$.)  Thus, for each $s\in (0,r)\setminus\{x_0\}$, there is exactly one solution $u$ such that $u(s)=0$ and $u(x_0)=1$.  Letting $s\to r$, the derivative $u'(0)$ converges to a negative limit, yielding a solution which is positive and decreasing on $(x_0,r)$.

If $v$ is any solution to $\mL v=\ga v$,
\begin{equation}  \label{E:integral}
v'(x)=e^{-B(x)}\left\{ v'(x_0) +2\int_{x_0}^x\lp \ga+\kappa(z)\rp v(z) dz\right\}.
\end{equation}
A consequence is that $v'(x)e^{B(x)}$ always has a limit at both boundaries.

Let us take a solution with $v(x_0)=0$ and $v'(x_0)>0$.  All terms in \eqref{E:integral} are positive for $x>x_0$, so $v$ is nondecreasing.  It follows that
$$
0\le v'(x)\le v'(x_0)e^{-B(x)} + 2(\ga+K)v(x)e^{-B(x)} \int_{x_0}^x e^{B(z)} dz=: B_1(x)+v(x)B_2(x).
$$

Suppose $r$ is regular.  Since $0$ and $r$ are regular, $B_1$ and $B_2$ are integrable on $(0,r)$.  It follows that $v$ is bounded on $(0,r)$.  From this $v$, and the nonincreasing $u$ mentioned earlier, we can compose two positive nonincreasing solutions $u_0$ and $u_1$ on all of $(0,r)$, such that
\begin{xalignat*}{2}
        \lim_{x\to r} u_0(x)&=0, &      \lim_{x\to r} e^{B(x)}u'_0(x)&=-1,\\
        \lim_{x\to r} u_1(x)&=1 &       \text{ and }\lim_{x\to r} e^{B(x)}u'_1(x)&=0.\\
\end{xalignat*}
Both solutions will be positive at $0$.  As has already been pointed out, the functions $v_i:=e^B u_i$ are eigenfunctions for the adjoint operator, satisfying $\mL^* v_i=\ga v_i$.  Since they are linearly independent, they span the solution space, and we can represent $\phi$ as $a_0 v_0+a_1 v_1$ for some real numbers $a_0$ and $a_1$.  By \eqref{E:fellerboundr}, if $a_0$ and $a_1$ are not both zero,
$$
0\ge\lim_{x\to r}  \frac{e^{-B(x)}\phi(x)}{ \oh\phi'(x)-b(x)\phi(x)}\frac{a_1}{-a_0}.
$$
This implies that $a_1$ and $a_0$ have the same sign, or $a_1=0$.  Without loss of generality, suppose both coefficients are positive.  We then have
$$
0\le \lim_{x\to 0} \frac{e^{-B(x)}\phi(x)}{ \oh\phi'(x)-b(x)\phi(x)}= \frac{a_0 u_0(0)+a_1 u_1(0)}{a_0 e^{B(x)}u'_0(0)+a_1 e^{B(x)}u'_1(0)}.
$$
Here the numerator is positive, and the denominator negative, which is a contradiction, proving that $a_0$ and $a_1$ are indeed both zero.

Consider now the situation when $r=\infty$ is a natural boundary.  Let $u$ be an eigenfunction as above, with $u(x_0)\ge 0$ and $u'(x_0)>0$.  By \eqref{E:integral},
$$
u'(x)\ge \ga e^{-B(x)} \int_{x_0}^x e^{B(z)} dz.
$$
As $\infty$ is a natural boundary, it must be that $e^{-B(x)}\int_{x_0}^x e^{B(z)}dz$ is not integrable. Together with the fact that $u$ is increasing, this implies that $\lim_{x\to\infty}u(x)=\infty$. Furthermore, for $x\ge x_1>x_0$ we have $\beta>0$ such that
$$
u(x)\ge \beta \int_{x_1}^x e^{-B(z)} dz.
$$
Again, since $\infty$ is a natural boundary, it follows that
$$
\int^\infty  e^{B(x)} u(x) dx =\infty.
$$
Since $\phi\in L^1$, it must be that $u=e^{-B}\phi$ is positive and nonincreasing, negative and nondecreasing, or identically 0.  By the same computation as above, either of the first two cases implies that
$$
\lim_{x\to 0} \frac{e^{-B(x)}\phi(x)}{ \phi'(x)-2b(x)\phi(x)}\le 0,
$$
contradicting the assumed boundary condition \eqref{E:fellerbound}.  Thus $\phi\equiv 0$.

Finally, suppose $r$ is an entrance boundary.  By the same argument as above, we see that solutions $u$ passing through 0 are unbounded at $r$, and when $u'(x_0)>0$ for $u(x_0)=0$, we have $\lim_{x\to r} e^{B(x)}u'(x)>0$.  It follows that, for $\phi$ to satisfy the boundary condition \eqref{E:fellerboundentrance} it must be a multiple of $e^B u_0$, where $u$ is a nonincreasing solution.  As before, though, a nonincreasing solution cannot satisfy the boundary condition \eqref{E:fellerbound} at 0.
\end{proof}

\subsection{Facts used in the proof of Theorem \ref{T:limsups}} \label{sec:tls}
\begin{Lem}
\label{L:ratiowork}
Suppose (LP) holds.  Then for any distribution $\nu$ satisfying (IDC), and any positive $z$ and $s$,
    \begin{equation}  \label{L:limitdecline}
        \lim_{t\to\infty}\frac{\P_{\nu}\{ X_{t+s}\le z \}}{\P_{\nu} \{X_{t}\le z \}} = e^{-\ulm s}. 
        \end{equation}
The limit is uniform in $s$.
\end{Lem}

\begin{proof}
Using the notation of section \ref{sec:mainnatreg}, we need to show that
\begin{equation}  \label{E:lhsrhs}
\ltf \frac{\mI(t+s,g\xi_{[0,z]})}{\mI(t,g\xi_{[0,z]})}=e^{-\ulm s},
\end{equation}
uniformly in $s$.  Let 
$$
h(\gl)=\sup_{s} \lp e^{-\ulm s} - e^{-\gl s}\rp\le \frac{\gl}{\ulm}-1.
$$
For all $t$ large enough that $\mI(t,g\xi_{[0,z]})>0$,
\begin{align*}
 \sup_{s\in [0,\infty)}\left|\frac{\mI(t+s,g\xi_{[0,z]})}{\mI(t,g\xi_{[0,z]})} 
- e^{-\ulm s} \right|&\sup_{s\in [0,\infty)}\left|\frac{\mI(t+s,g\xi_{[0,z]})}{\mI(t,g\xi_{[0,z]})} 
- \frac{\mI(t,g\xi_{[0,z]}e^{-\ulm s})}{\mI(t,g\xi_{[0,z]})} \right|\\
&\le \frac{\mI(t,hg\xi_{[0,z]})}{\mI(t,g\xi_{[0,z]})}
\end{align*}
An application of \eqref{E:uselater} proves that this converges to 0 as $t\to\infty$.
\end{proof}

\begin{Lem}  \label{L:twolims}
    Let $f,c,\gm:\R^{+}\times \Z^{+}\to\R^{+}$ be maps satisfying
    \begin{equation}
        f(z,n+1)\le c(z,n)f(z,n)+\gm(z,n),
        \label{E:fineq}
    \end{equation}
    and $c:=\lzn c(z,n)<1$. Let $\gm:=\lzn \gm(z,n)$ Then
    \begin{equation}
        \lzn f(z,n)\leq\frac{\gm}{(1-c)}.
        \label{E:fbound}
    \end{equation}
\end{Lem}

\begin{proof}
    Choose $\gep\in(0,1-c)$, and find $n_{0}$ and $z_{0}$ such that 
    for any $z\ge z_{0}$ and $n\ge n_{0}$, $c(z,n)\le c+\gep$ and 
    $\gm(z,n)\le \gm+\gep$.  Then for any $z\ge z_{0}$ and $n\ge n_{0}$ 
    repeated application of \eqref{E:fineq} yields
    \begin{align*}
    f(z,n)\le (c+\gep)^{n-n_{0}} &f(z,n_{0})+(\gm+\gep)\lp 
    1+(c+\gep)+\cdots+(c+\gep)^{n-n_{0}-1} \rp\\
    &=(c+\gep)^{n-n_{0}} 
    f(z,n)+(\gm+\gep)\frac{1-(c+\gep)^{n-n_{0}}}{1-c-\gep}.
    \end{align*}
    Thus
    $$
    \lzn f(z,n) \le \frac{\gm+\gep}{1-c-\gep}.
    $$
    Since $\gep$ was arbitrary, this completes the proof.
\end{proof}

\begin{Lem}
\label{L:subinvar}
If $\izf \plm(y)dy$ is finite, then $\plm$ is $\ulm$-subinvariant for the semigroup $P^{*}_{s}$.  That is,
\begin{equation}  \label{E:wrongway}
P^{*}_{s}(\plm)\le\plm P^{*\phi_{\ulm}}_{s} \indic \le e^{-\ulm s}\phi.
\end{equation}
If
\begin{equation}  \label{E:superinv}
\izf \P_{y}\ls\tp >s\rs\phi_{\ulm}(y)dy \ge e^{-\ulm s},
\end{equation}
then $\plm$ is $\ulm$-invariant.
\end{Lem}

\begin{proof}
\nc{\htr}{\mL^{*\phi_{\ulm}}}
Define the h-transformed operator
$$
\htr f = \rec{\plm} \mL^{*} (f\plm) = 
\oh \frac{d^{2}f}{dx^{2}} + \Bigl( \frac{\plm'}{\plm}-b\Bigr)\frac{df}{dx} -\ulm f.
$$
This defines a diffusion on $\R^{+}$ with constant killing at rate $\ulm$.  If we let $P_{s}^{*\phi_{\ulm}}$ be the corresponding semigroup, we have 
$$
P^{*}_{s}\plm\le\plm P^{*\phi_{\ulm}}_{s} \indic \le e^{-\ulm s}\phi.
$$
(This is an inequality rather than equality because the diffusion could explode.)  If \eqref{E:superinv} holds as well, then
$$
\izf \lb P^{*}_{s}\plm(y) - e^{-\ulm s}\phi(y) \rb dy\ge 0 \ge \izf P^{*}_{s}\plm(y)dy - e^{-\ulm s}.
$$
Since the right-hand side and left-hand side are equal, it must be that the integrand is zero almost everywhere.  We conclude by continuity that $P^{*}_{s}\plm(y)=e^{-\ulm s}\phi(y)$ for all $y$.
\end{proof}

Note that this is essentially a special case of Theorem 4.8.5 of \cite{rP95}.

\section*{Acknowledgement}
The authors would like to thank Persi Diaconis, Jim Pitman, and Jaime San Mart\'in, who provided helpful comments and guidance to the literature.

\sng

\end{document}